\documentclass[times,a4paper]{article}

\usepackage[separate-uncertainty=true,multi-part-units=single,allow-number-unit-breaks]{siunitx}

\usepackage[margin=3cm]{geometry}
\usepackage{graphicx}
\usepackage{authblk}
\usepackage{upgreek}
\usepackage{textcomp}
\usepackage{stmaryrd}
\usepackage{verbatim} 
\usepackage{float, booktabs, bm, pgfplots, subcaption}
\usepackage{amsmath, amsfonts, amssymb, times}
\usepackage[vlined,linesnumbered]{algorithm2e}
\usepackage{todonotes}
\usepackage[capitalize]{cleveref}
\usepackage{multirow}

\usepgfplotslibrary{colorbrewer}

\pgfplotsset{
  log x ticks with fixed point/.style={
      xticklabel={
        \pgfkeys{/pgf/fpu=true}
        \pgfmathparse{exp(\tick)}
        \pgfkeys{/pgf/fpu=false}
      }
  },
  log y ticks with fixed point/.style={
      yticklabel={
        \pgfkeys{/pgf/fpu=true}
        \pgfmathparse{exp(\tick)}%
        \pgfmathprintnumber[fixed relative, precision=3]{\pgfmathresult}
        \pgfkeys{/pgf/fpu=false}
      }
  }
}

\pgfplotsset{width=8.0cm,
	     scaled y ticks=false,
             cycle list/Dark2,
             cycle multiindex* list={
	       mark list*\nextlist
	       linestyles\nextlist
	       Dark2\nextlist
	     },
	    }

\pgfplotsset{
  log x ticks with fixed point/.style={
      xticklabel={
        \pgfkeys{/pgf/fpu=true}
        \pgfmathparse{exp(\tick)}%
        \pgfmathprintnumber[fixed relative, precision=3]{\pgfmathresult}
        \pgfkeys{/pgf/fpu=false}
      }
  },
  log y ticks with fixed point/.style={
      yticklabel={
        \pgfkeys{/pgf/fpu=true}
        \pgfmathparse{exp(\tick)}%
        \pgfmathprintnumber[fixed relative, precision=3]{\pgfmathresult}
        \pgfkeys{/pgf/fpu=false}
      }
  }
}

\newcommand{\T}{\mathrm{T}}
\newcommand{\e}{\mathrm{e}}

\newcommand{\p}{\mathrm{p}}

\newcommand{\comp}{\mathrm{c}}
\newcommand{\tens}{\mathrm{t}}
\newcommand{\jump}{\llbracket u \rrbracket}
\newcommand{\eq}{\mathrm{eq}}

\newcommand{\tr}{\mathrm{tr}}
\newcommand{\etal}{\textit{et al}. }
\newcommand{\ie}{\textit{i}.\textit{e}. }
\newcommand{\eg}{\textit{e}.\textit{g}. }

\newcommand{\real}{\mathbb{R}}

\newcommand{\Ful}{\ensuremath{\mathcal{A}}}
\newcommand{\added}{\ensuremath{\mathcal{T}}}
\newcommand{\ful}{\mathrm{f}}

\newcommand{\norm}[1]{\left\lVert{#1}\right\rVert}

\newcommand{\mrm}{\mathrm}
\newcommand{\mbf}{\mathbf}
\newcommand{\bs}{\boldsymbol}

\newcommand{\prob}[1]{\ensuremath{p\left(#1\right)}}

\newcommand{\fetwo}{FE$^\mathrm{2}$}
\newcommand{\strain}{\ensuremath{\bs{\varepsilon}}}
\newcommand{\stress}{\ensuremath{\bs{\sigma}}}
\newcommand{\disp}{\ensuremath{\mbf{u}^\omega}}

\newcommand{\basis}{\ensuremath{\boldsymbol{\Phi}}}

\newcommand{\dataset}{\ensuremath{\mathcal{D}}}
\newcommand{\weights}{\ensuremath{\mathbf{w}}}
\newcommand{\xvec}{\ensuremath{\mathbf{x}}}
\newcommand{\Xvec}{\ensuremath{\mathbf{X}}}
\newcommand{\yvec}{\ensuremath{\mathbf{y}}}
\newcommand{\noise}{\ensuremath{\sigma^2}}
\newcommand{\gtol}[1]{\ensuremath{\gamma_\mathrm{tol}=\SI{#1}{\mega\pascal}}}
\newcommand{\sigf}[1]{\ensuremath{\sigma_\mathrm{f}^2=\SI{#1}{\mega\pascal\squared}}}
\newcommand{\length}[1]{\ensuremath{\ell=\SI{#1}{}}}
\newcommand{\sign}[1]{\ensuremath{\sigma_\mathrm{n}^2=\SI{#1}{\mega\pascal\squared}}}
\newcommand{\obs}{\ensuremath{\mathrm{o}}}

\begin{document}

\title{On-the-fly construction of surrogate constitutive models for concurrent multiscale mechanical analysis through probabilistic machine learning}

\author[1]{I. B. C. M. Rocha}
\author[2,3]{P. Kerfriden}
\author[1]{F. P. van der Meer}

\affil[1]{Delft University of Technology, Faculty of Civil Engineering and Geosciences, P.O. Box 5048, 2600GA Delft, The Netherlands}
\affil[2]{Mines ParisTech (PSL University), Centre des mat\'{e}riaux, 63-65 Rue Henri-Auguste Desbru\`{e}res BP87, F-91003 \'{E}vry, France}
\affil[3]{Cardiff University, School of Engineering, Queen's Buildings, The Parade, Cardiff, CF24 3AA, United Kingdom}

\date{}
\maketitle

\begin{abstract}
Concurrent multiscale finite element analysis (\fetwo) is a powerful approach for high-fidelity modeling of materials for which a suitable macroscopic constitutive model is not available. However, the extreme computational effort associated with computing a nested micromodel at every macroscopic integration point makes \fetwo\ prohibitive for most practical applications. Constructing surrogate models able to efficiently compute the microscopic constitutive response is therefore a promising approach in enabling concurrent multiscale modeling. This work presents a reduction framework for adaptively constructing surrogate models for \fetwo\ based on statistical learning. The nested micromodels are replaced by a machine learning surrogate model based on Gaussian Processes (GP). The need for \emph{offline} data collection is bypassed by training the GP models \emph{online} based on data coming from a small set of fully-solved \emph{anchor} micromodels that undergo the same strain history as their associated macroscopic integration points. The Bayesian formalism inherent to GP models provides a natural tool for \emph{online} uncertainty estimation through which new observations or inclusion of new \emph{anchor} micromodels are triggered. The surrogate constitutive manifold is constructed with as few micromechanical evaluations as possible by enhancing the GP models with gradient information and the solution scheme is made robust through a greedy data selection approach embedded within the conventional finite element solution loop for nonlinear analysis. The sensitivity to model parameters is studied with a tapered bar example with plasticity, while the applicability of the model to more complex cases is demonstrated with the elastoplastic analysis of a plate with multiple cutouts and a crack growth example for mixed-mode bending. The framework is found to be a promising approach in reducing the computational cost of \fetwo, with significant efficiency gains being obtained without resorting to \emph{offline} training.
\end{abstract}

\textbf{Keywords:}
Concurrent multiscale, Surrogate modeling, Probabilistic learning, Gaussian Processes (GP), Active learning.

\section{Introduction}
\label{SEintroduction}

There is a growing demand for high-fidelity numerical techniques capable of describing material behavior across spatial scales. With recent advances in additive manufacturing allowing for the development of novel materials with highly-tailored microstructures \cite{Gantenbein2018}, multiscale modeling techniques will become increasingly relevant in the design of novel materials and structures. One popular approach for concurrent multiscale modeling is the so-called \fetwo\ approach \cite{Feyel1999,Geers2010}, in which macroscopic material response is directly upscaled from embedded micromodels without introducing additional constitutive assumptions. \fetwo\ is a powerful and versatile technique used in a number of solid mechanics applications, from continuous \cite{Kouznetsova2001} and discontinuous \cite{Nguyen2012} mechanical equilibrium to multiphysics problems involving heat and mass transfer \cite{Ozdemir2008} and material degradation due to aging \cite{Terada2010,Rocha2017a}. However, \fetwo\ has the major drawback of being associated with extreme computational costs, hindering its application in actual design scenarios. Enabling the use of \fetwo\ in many practical applications that would benefit from its accuracy and versatility is therefore highly contingent on being able to reduce its computational cost to tractable levels.

A promising approach in accelerating \fetwo\ models consists in constructing surrogate models that take the place of the original high-fidelity micromodels at each macroscopic integration point. When building surrogates, the goal is to maintain as much of the generality offered by the original micromodel while eliminating as much computational complexity as possible. One option is to employ unsupervised learning on a number of full-order solution snapshots in order to define lower-dimensional solution manifolds for both displacements \cite{Kunisch2002,Kerfriden2011} and internal forces \cite{Hernandez2017,vanTuijl2019} at the microscale \cite{Liu2016,Goury2016,Rocha2018b}. Alternatively, a supervised learning approach can be taken by using snapshots of the homogenized micromodel response to directly define a data-driven regression model for the macroscopic constitutive behavior \cite{Le2015,Lu2018,Ghavamian2019,Liu2019,Mozaffar2019,Rocha2020c}. Although resulting in models of distinct natures, both approaches rely on the existence of an observation database on the behavior of the original micromodel that is usually obtained \emph{offline} (before deployment on a multiscale setting) and should cover every possible scenario the surrogate is expected to approximate \emph{online}.

However, building such a database of model snapshots can be a challenging task (see \cite{Bessa2017,Mozaffar2019} for interesting approaches based on Design of Experiments and \cite{Goury2016} for a data-driven training framework based on Bayesian Optimization). For micromodels employing path-dependent materials, this \emph{offline training} process entails sampling a highly-complex constitutive manifold that depends on an arbitrarily long strain history and can therefore be excessively sensitive to small changes in boundary conditions (\eg strain localization and crack propagation problems). Furthermore, training a surrogate with such a complex dataset often requires additional partitioning techniques in order to avoid computationally inefficient reduced models \cite{Kerfriden2013,Ghavamian2017}. An alternative to the conventional \emph{offline-online} approach that has been gaining popularity is the use of adaptive reduction frameworks that either preclude the need for \emph{offline} training altogether \cite{Ryckelynck2005,Rocha2020a} or combine different reduction techniques into a single framework with only limited \emph{offline} effort while employing \emph{online} error indicators to continuously assess the quality of the approximation and trigger a refinement of the surrogates when necessary. Nevertheless, obtaining consistent adaptivity criteria for either hyper-reduced models \cite{Ryckelynck2005} or machine learning models based on least-squares solutions \cite{Fritzen2019,Liu2019} is not straightforward. At the other end of the spectrum, works dealing with constitutive models based on Bayesian regression techniques, that provide a natural way to estimate error \cite{Tartakovsky2018,Salloum2014,Frankel2019}, do not take advantage of their potential for creating adaptive frameworks. There is a need, therefore, for the development of fully-online approaches with reliable adaptivity strategies based on sequential learning techniques and on error estimation methods with robust probabilistic foundations.

\begin{figure}[htb]
\centering
\includegraphics[width=0.8\textwidth]{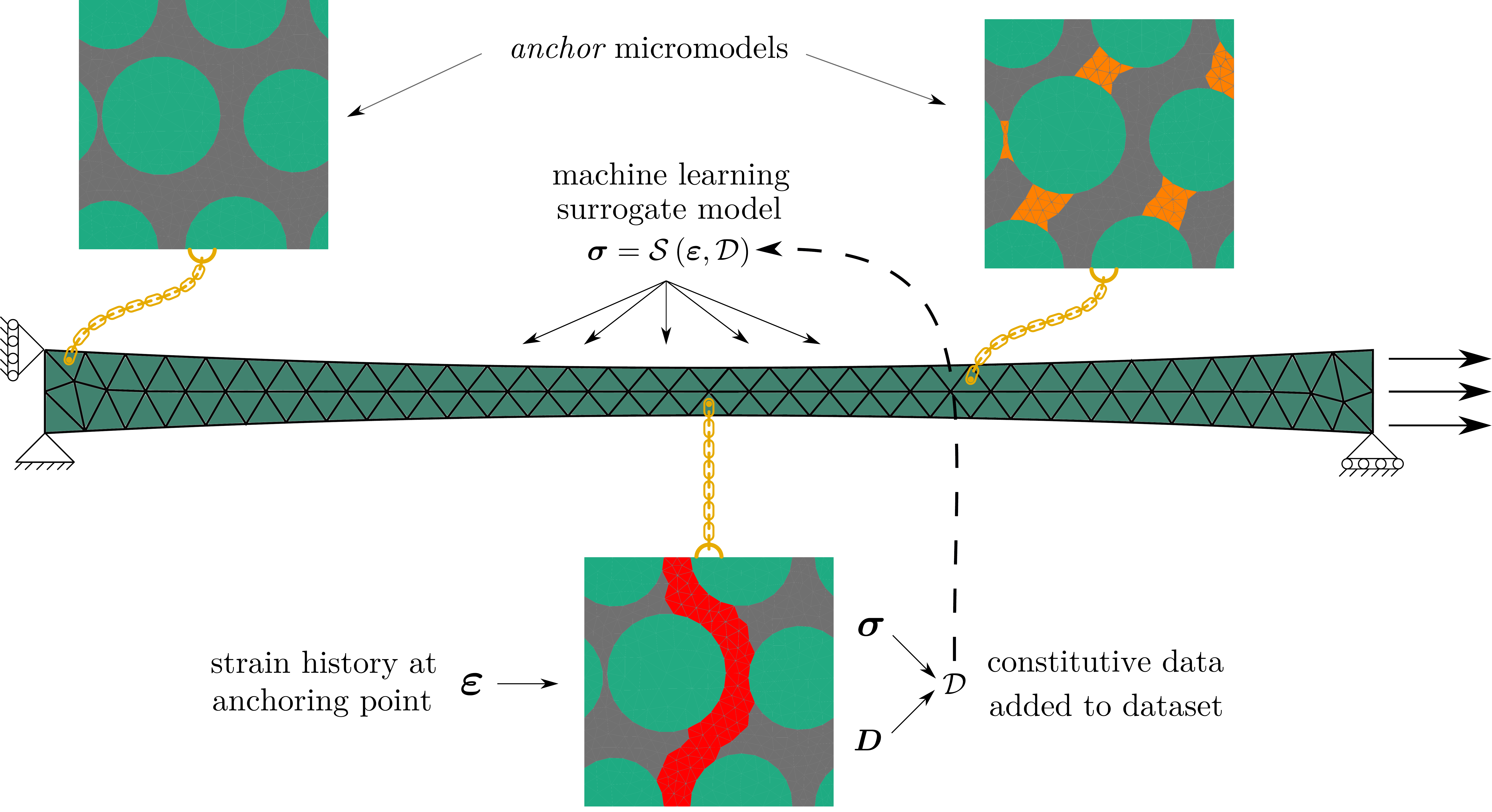}
\caption{Schematic representation of the \emph{online} adaptive reduction framework presented in this work. A small number of \emph{anchor} models are used to train a machine learning surrogate that evolves as the macroscopic structure is loaded.}
\label{FIintro}
\end{figure}

This work presents an adaptive probabilistic framework for constructing surrogate constitutive models for nonlinear concurrent multiscale analysis. The approach is based on substituting the original models associated to macroscopic integration points with a machine learning surrogate. In order to bypass the need for an \emph{offline} training phase, a small number of fully-solved models associated to representative macroscopic integration points are used to generate constitutive data \emph{online}. Because these models are not directly used to make predictions at every macroscopic iteration but only provide an indirect coupling between the scales, we denote them as \emph{anchor models} (\cref{FIintro}). Stress and stiffness data resulting from subjecting the \emph{anchor} models to the same strain histories seen by their respective anchoring points is used to build a surrogate model based on the Gaussian Process (GP) regression technique. The accuracy of this reduced-order solution is controlled with uncertainty information that arises naturally from the Bayesian formalism of the GP models. The resultant material model is embedded within a conventional finite element solution in a way that ensures the macroscopic solution is numerically robust and limits the sampling of new data as much as possible in order to maximize efficiency. The framework is demonstrated with a set of numerical tests including both bulk strain localization and crack propagation in order to assess its accuracy, efficiency and versatility.

\section{Concurrent multiscale analysis}
\label{SEconcurrentmultiscaleanalysis}

\subsection{Macroscopic problem}
\label{SEmacroscopicproblem}

We begin by briefly introducing the concurrent multiscale equilibrium problem that we seek to accelerate. Let $\Omega$ define the macroscopic domain being modeled. We wish to find the displacement field $\mbf{u}^\Omega$ resulting from a set of Dirichlet and Neumann boundary conditions applied to the surface $\Gamma$ that bounds $\Omega$. Under the assumption of small strains, the equilibrium solution is found by satisfying: 
\begin{equation}
\mrm{div}\left(\boldsymbol{\sigma}^\Omega\right) = \mbf{0}
\quad\quad
\strain^\Omega = \frac{1}{2}\left(\nabla\mbf{u}^\Omega+\left(\nabla\mbf{u}^\Omega\right)^\T\right)
\end{equation}
\noindent where $\mrm{div}\left(\cdot\right)$ is the divergence operator and body forces have been neglected. In order to solve for $\mbf{u}^\Omega$, a constitutive model $\mathcal{M}$ that relates $\bs\sigma^\Omega$ and $\bs\varepsilon^\Omega$ must be introduced:
\begin{equation}
\bs\sigma^\Omega = \mathcal{M}\left(\bs\varepsilon^\Omega,\bs\varepsilon^\Omega_\mrm{h}\right)
\label{EQfe2macroconstmodel}
\end{equation}
\noindent where $\bs\varepsilon^\Omega_\mrm{h}$ is a history term that accounts for strain path dependency. In the context of multiscale analysis, the model $\mathcal{M}$ can be seen as a homogenization operator that lumps all physical processes happening at scales lower than $\Omega$ into a homogeneous medium with equivalent behavior. Depending on how complex the microscopic behavior is, $\mathcal{M}$ can range from having a relatively simple form (\eg linear elasticity) to being next to impossible to formulate explicitly. 

\subsection{Microscopic problem}
\label{SEmicroscopicproblem}

In a concurrent multiscale approach, we do not formulate $\mathcal{M}$ directly and instead opt for upscaling microscopic behavior from micromodels embedded at each macroscopic material point. Let $\omega$ be a Representative Volume Element (RVE) of the microscopic material features whose behavior we would like to upscale. Assuming the principle of separation of scales is valid ($\omega\ll\Omega$), we can link the two scales by enforcing:
\begin{equation}
\disp = \strain^\Omega\mbf{x}^\omega + \widetilde{\mbf{u}}
\label{EQfe2disp}
\end{equation}
\noindent where the response is decomposed into a linear displacement field imposed by the presence of macroscopic strains and a fluctuation field $\widetilde{\mbf{u}}$ that accounts for microscopic inhomogeneities. 

In order to obtain $\bs\sigma^\Omega$, we first find an equilibrium solution for $\mbf{u}^\omega$ by satisfying:
\begin{equation}
\mrm{div}\left(\stress^\omega\right)=\mbf{0}
\quad\quad
\strain^\omega = \frac{1}{2}\left(\nabla\mbf{u}^\omega+\left(\nabla\mbf{u}^\omega\right)^\T\right)
\end{equation}
\noindent where we see that once more the need arises for the definition of a constitutive model, this time relating $\bs\varepsilon^\omega$ with $\bs\sigma^\omega$. Underlying the choice for a multiscale approach is the assumption that constitutive behavior can be represented by models of decreasing complexity as one descends to lower scales. It is therefore common to employ regular constitutive models (\eg (visco)elasticity, (visco)plasticity, damage) for material components at the microscale \cite{Kouznetsova2001,Ozdemir2008,Terada2010,Rocha2017a}. However, the framework is flexible in allowing for models on an even lower third scale to be embedded to material points in $\omega$ (at the cost of even higher computational effort).

\subsection{Bulk homogenization}
\label{SEbulkhomogenization}

With a solution for $\mbf{u}^\omega$ under the constraint that the fluctuation field $\widetilde{\mbf{u}}$ must be periodic, we can use the Hill-Mandel principle \cite{Hill1972} to obtain homogenization expressions for the macroscopic stresses $\bs\sigma^\Omega$ and consistent tangent stiffness $\mbf{D}^\Omega$:
\begin{equation}
\stress^\Omega = \frac{1}{\omega}\int_\omega\stress^\omega\mrm{d}\omega
\quad\quad
\mbf{D}^\Omega = \mathcal{P}\left(\mbf{K}^\omega\right)
\label{EQfe2stress}
\end{equation}
\noindent from which we obtain the intuitive result that the macroscopic stresses are simply the volume average of the microscopic ones. Finally, the macroscopic constitutive tangent stiffness is computed through a probing operator $\mathcal{P}$ applied on the global microscopic tangent stiffness matrix $\mbf{K}^\omega$ \cite{Nguyen2012}.

\subsection{Cohesive homogenization}
\label{SEcohesivehomogenization}

Although the preceding formulation allows for general microscopic constitutive behavior to be upscaled, the response loses objectivity with respect to the RVE size after the onset of global microscopic softening \cite{Nguyen2010}, \ie when the determinant of the acoustic tensor is zero along a given direction $\mbf{n}$:
\begin{equation}
\det\left(\mbf{n}^\T\mbf{D}^\Omega\mbf{n}\right)=0
\label{EQacoustictensor}
\end{equation}
\noindent This non-objectivity arises because the volume $\omega_\mrm{d}$ of the strain localization band that causes the softening does not scale with the RVE size, an observation that motivates the use of a modified version of the Hill-Mandel principle \cite{Nguyen2010,Nguyen2012}:
\begin{equation}
\frac{1}{w}\bs\tau^\Omega\delta\mbf{\llbracket v \rrbracket}^\Omega
=
\frac{1}{wh}\int_{\omega_\mrm{d}}\bs\sigma^\omega\delta\bs\varepsilon^\omega\mrm{d}\omega
\label{EQfe2traction}
\end{equation}
\noindent where $w$ and $h$ are geometric RVE parameters that depend on the localization band orientation, $\bs\tau^\Omega$ is a macroscopic traction and $\mbf{\llbracket v \rrbracket}^\Omega$ is a shifted displacement jump that allows for an initially-rigid cohesive response. Note that the homogenization is now performed towards a cohesive traction acting on a macroscopic surface that defines a discontinuity in $\mbf{u}^\Omega$.

The development of a consistent strategy for continuous-discontinuous scale linking in \fetwo\ is an open issue that is left out of the scope of the present discussion, with a number of different approaches being found, for instance, in \cite{Nguyen2010,Nguyen2012,Sanchez2013,Oliver2015,Svenning2019}. For the purpose of building surrogate models for the RVE response, it suffices to acknowledge that two distinct models should be trained for bulk and cohesive responses, as the underlying constitutive manifolds have dimensions with different physical interpretations (strain/stress versus jump/traction).

\subsection{Acceleration strategy}
\label{SEaccelerationstrategy}

Assuming FEM is used to solve the mechanical problems at both scales, we can approximate the computational cost of a single macroscopic iteration as:
\begin{equation}
\mrm{Cost} \approx \left(N_\mrm{dof}^\Omega\right)^{x}+N_\mrm{ip}^\Omega N_\mrm{iter}^\omega\left(N^\omega_\mrm{dof}\right)^{x}
\label{EQfe2cost}
\end{equation}
\noindent where $N_\mrm{ip}^\Omega$ is the number of macroscopic integration points and $N_\mrm{iter}^\omega$ is the number of iterations necessary for convergence of the microscopic BVP. We assume for simplicity that the bulk of the effort comes from solving the linearized systems of equations involving $\mbf{K}^\Omega$ and $\mbf{K}^\omega$ for nodal displacements, where the complexity exponent $x$ depends on the solver used. It can be seen that the second term, associated with solving the microscopic equilibrium problems, quickly outweighs the first and becomes a performance bottleneck as $N_\mrm{dof}^\omega$ increases, especially since the number of macroscopic integration points $N_\mrm{ip}^\Omega$ increases together with $N_\mrm{dof}^\Omega$. Constructing a surrogate that replaces the original constitutive model $\mathcal{M}$ of \cref{EQfe2macroconstmodel} is therefore an effective approach to accelerating \fetwo.

However, constructing such a surrogate \emph{offline} is a challenging undertaking, otherwise there would not have been need for a multiscale approach in the first place. From \cref{EQfe2macroconstmodel} we see that the constitutive manifold to be reproduced can have an arbitrarily high dimensionality due to the dependency on $\bs\varepsilon^\Omega_\mrm{h}$. This is equivalent to stating that the shape of the $\bs\varepsilon^\Omega$-$\bs\sigma^\Omega$ manifold can change after each load step. Sampling this high-dimensional input space \emph{offline} in order to have a surrogate that is accurate for arbitrary strain histories seems to be an intractable problem that, to the best of our knowledge, has still not been tackled in a satisfactory way. 

This issue can be avoided by exploiting the fact that $\bs\varepsilon^\Omega$ (and therefore $\bs\varepsilon^\Omega_\mrm{h}$) are often highly constrained by the geometry and boundary conditions of the macroscopic structure being modeled \cite{Ghavamian2017}. We therefore opt for constructing a highly-tailored surrogate model $\mathcal{S}$ \emph{online} based on a dataset \dataset\ of observations coming from a small number of fully-solved micromodels:
\begin{equation}
\bs\sigma^\Omega = \mathcal{S}\left(\bs{\varepsilon}^\Omega,\dataset\right)
\label{EQsurrogateconstmodel}
\end{equation}
\noindent which can then be used to compute the constitutive response for a fraction of the cost. When trying to keep \dataset\ as small as possible for the case at hand, it is crucial to have a means for quantifying the uncertainty in probing $\mathcal{S}$ for any given $\strain^\Omega$. In this work, a Bayesian approach is adopted to assess \emph{online} whether \dataset\ is large enough to provide the desired level of confidence in $\mathcal{S}$ at a given $\strain^\Omega$. 

\section{Bayesian surrogate modeling}

In this section we introduce the Bayesian regression approach used to construct surrogate constitutive models, beginning from parametric versions of the surrogate model $\mathcal{S}$ --- \ie by encapsulating the constitutive information in \dataset\ into a set of parameters $\mbf{w}$ --- and eventually moving to a non-parametric model based on Gaussian Processes (GP) that uses the data in \dataset\ directly in order to make predictions. Our goal here is to appeal to the reader who might be unfamiliar with probabilistic regression models by starting from classical least-squares regression and gradually moving towards a Bayesian approach. Nevertheless, the discussion is kept as brief and focused as possible. The interested reader can find richer discussions on the subject in \cite{Bishop,Rasmussen}.

\subsection{Least-squares regression}

We start by building a parametric model $y\left(\weights,\xvec\right)$ that approximates a scalar target response $t$ by fitting \weights\ with a dataset \dataset\ of $N$ observations $\mbf{t}_\obs$ at $\mbf{X}_\obs=\left[\mbf{x}_{\obs 1} \cdots \mbf{x}_{\obs N}\right]$. We therefore assume that the target $t$ can be written as:
\begin{equation}
t = y\left(\weights,\xvec\right) + \epsilon
\quad\mrm{with}\quad
y\left(\weights,\xvec\right)=\sum_j^M w_j\phi_j\left(\xvec\right)
\,\,\mrm{and}\,\,
\prob{\epsilon}=\mathcal{N}\left(\epsilon\vert0,\noise_\mrm{n}\right)
\label{EQnoisylinregression}
\end{equation}
\noindent where the noise $\epsilon$ is given by a zero-mean Gaussian distribution with variance $\noise_\mrm{n}$ and $y$ is linear with respect to its $M$ weights, gathered in the vector \weights. Note that the assumption of a linear model does not limit its fitting capabilities since the basis functions $\bs\phi$ define a feature space that can be nonlinear in $\mbf{x}$. In the context of this work, $t$ can be a single stress or traction component, but the discussion is equally applicable to the problem of modeling individual model components (\eg yield parameters). 

In order to find values for \weights, we compute the \emph{likelihood function} of the model, \ie how likely the model is to produce the values $\mbf{t}_\obs$ in \dataset\ given \weights: 
\begin{equation}
\prob{\mbf{t}_\obs\vert\weights,\noise_\mrm{n}}=\prod_i^N\mathcal{N}\left({t}_{\obs i}\vert\weights^\T\bs{\phi}\left(\xvec_{\obs i}\right),\noise_\mrm{n}\right)
\label{EQlikelihood}
\end{equation}
\noindent which is a product of the probabilities of each point in isolation because we assume that each sample $t_{\obs i}$ is sampled from the conditional distribution $\prob{t\vert y}$ independently. We can find an optimum data fit for \weights\ by maximizing the likelihood and assuming that the shapes of $\bs\phi$ are fixed \emph{a priori}: 
\begin{equation}
\nabla_\mbf{w}\prob{\mbf{t}_\obs\vert\weights,\noise_\mrm{n}}=0 \quad\Rightarrow\quad \weights_\mrm{ML}=\left(\basis^\T\basis\right)^{-1}\basis^\T\mbf{t}_\obs\,\,\mrm{and}\,\,\noise_\mrm{n}=\frac{1}{N}\sum_i^N\left(t_{\obs i}-\weights_\mrm{ML}^\T\bs{\phi}\left(\xvec_{\obs i}\right)\right)^2
\label{EQmlequations}
\end{equation}
\noindent where $\basis\in\real^{N\times M}$ is a matrix with basis function values $\phi_{i,j}=\phi_i\left(\mbf{x}_{\obs j}\right)$ evaluated at each point in \dataset. Note that this is equivalent to minimizing the sum of squared differences between $y$ and $t$, so $\weights_\mrm{ML}$ is the same parameter vector obtained by the classical least squares approach and $\noise_\mrm{n}$ quantifies the spread of $t$ around $y$.

Here we do not opt for a least-squares approach for three reasons. Firstly, it can suffer from severe overfitting when the dataset \dataset\ is small --- which is the case in the present work since we have no \emph{offline} training and only a small number of fully-solved micromodels to sample from. Secondly, the uncertainty associated with $\noise_\mrm{n}$ is constant throughout the input space $\xvec$ and does not provide an indication that the model is being used at a location far from data points (which could then be used to trigger a refinement of \dataset). Finally, performing model selection in the least-squares framework, using cross-validation for instance, is rather tedious compared to performing model selection in the context of Bayesian regression methods \cite{Bishop}.

\subsection{Bayesian parametric regression}
\label{SEbayesianparametricregression}

In the Bayesian approach to regression, we not only assume an uncertainty over the target $t$ but also over the weights \weights. We initially assume a prior probability over \weights\ that represents our initial model assumptions before any data is encountered:
\begin{equation}
\prob{\weights}=\mathcal{N}\left(\weights\vert\mbf{0},\noise_\mrm{w}\mbf{I}\right)
\label{EQweightprior}
\end{equation}
\noindent where $\noise_\mrm{w}$ is the variance parameter associated with the uncertainty over values of \weights. Information from \dataset\ is incorporated by using Bayes' theorem to obtain a posterior probability distribution for \weights:
\begin{equation}
\prob{\weights\vert\mbf{t}_\obs} = \frac{\prob{\mbf{t}_\obs\vert\weights}\prob{\weights}}{\prob{\mbf{t}_\obs}}
\quad\mrm{with}\quad
\prob{\mbf{t}_\obs}=\int\prob{\mbf{t}_\obs\vert\weights}\prob{\weights}\mrm{d}\weights
\label{EQbayestheorem}
\end{equation}
\noindent where $\prob{\mbf{t}_\obs\vert\weights}$ is the likelihood function of \cref{EQlikelihood} and $\prob{\mbf{t}_\obs}$ is the \emph{marginal likelihood} of the model (\ie the probability of producing the dataset \dataset). Because both \weights\ and $t\vert\weights$ are Gaussian variables, an analytical solution exists for $\prob{\weights\vert\mbf{t}_\obs}$ \cite{Bishop}:
\begin{equation}
\prob{\weights\vert\mbf{t}_\obs}=\mathcal{N}\left(\weights\vert\weights_\mrm{N},\mbf{S}_\mrm{N}\right)
\quad\mrm{with}\quad
\weights_\mrm{N}=\frac{1}{\noise_\mrm{n}}\mbf{S}_\mrm{N}\basis^\T\mbf{t}_\obs
\,\,\mrm{and}\,\,
\mbf{S}_\mrm{N}= \left(\frac{1}{\noise_\mrm{w}}\mbf{I}+\frac{1}{\noise_\mrm{n}}\basis^\T\basis\right)^{-1}
\label{EQmapequations}
\end{equation}
\noindent Note that the expected value of \weights\ now depends on both the data points in \dataset\ and on our initial beliefs about \weights\ represented by the prior distribution. Given this posterior, the best guess for \weights\ is the one with the highest $\prob{\weights\vert\mbf{t}_\obs}$. This is the so-called Maximum A Posteriori (MAP) value and, for a Gaussian posterior, $\weights_\mrm{MAP}=\weights_\mrm{N}$. This estimation for \weights\ is equivalent to a least-squares prediction with a quadratic regularization term proportional to $\noise_\mrm{w}$ which helps reducing the negative effects of overfitting.

However, in a fully Bayesian treatment of linear regression, we do not choose one specific value for \weights\ but rather make a prediction for a new target value $t_*$ by averaging over all possible values of \weights:
\begin{equation}
\prob{t_*\vert\mbf{t}_\obs,\noise_\mrm{n},\noise_\mrm{w}} = \int\prob{t_*\vert\weights,\noise_\mrm{n}}\prob{\weights\vert\mbf{t}_\obs,\noise_\mrm{w},\noise_\mrm{n}}\mrm{d}\weights
\label{EQlinmodelpredictivedistribution}
\end{equation}
\noindent 
In order to set the stage for Gaussian Processes, it is also interesting to represent the expectation of $t_*$ as:
\begin{equation}
\mathbb{E}\left[t_*\vert\xvec_*\right]=\sum_i^N k\left(\xvec_{\obs i},\xvec_*\right)t_{\obs i}
\quad\mrm{with}\quad
k\left(\xvec_p,\xvec_q\right)=\frac{1}{\noise_\mrm{n}}\bs{\phi}^\T\left(\xvec_p\right)\mbf{S}_\mrm{N}\bs{\phi}\left(\xvec_q\right)
\label{EQlinmodelkernelform}
\end{equation}
\noindent where $\weights$ has now vanished and the new prediction is cast as a linear combination of values from \dataset, where $k\left(\xvec_p,\xvec_q\right)$ is a \emph{kernel} that measures the similarity between two points in input space.

The Bayesian approach to the (generalised) linear model circumvents the problem of overfitting and provides confidence intervals. However, the linear model has the drawback of having a fixed number of basis functions $\bs\phi$ whose shapes need to be chosen before \dataset\ is observed. This can be circumvented by either allowing $\bs\phi$ to change in shape during training (\eg in a neural network, where $\bs\phi$ are neuron values coming from the last hidden layer) or by explicitly choosing a kernel $k\left(\xvec_p,\xvec_q\right)$ and adopting a distribution over functions instead of over weights, as in Gaussian Process regression models.

\subsection{Gaussian Process (GP) regression}
\label{SEgaussianprocesses}

Since the regression function $y$ is a function of \weights, adopting a prior distribution over \weights\ (\cref{EQweightprior}) implicitly leads to a distribution over function values $\mbf{y}$. This is the basis of Gaussian Process (GP) models. Here we assume these function values are jointly Gaussian:
\begin{equation}
\prob{\mbf{y}}=\mathcal{N}\left(\mbf{y}\vert\mbf{0},K\left(\mbf{X},\mbf{X}\right)\right)
\label{EQyprior}
\end{equation}
\noindent where a zero mean is assumed without loss of generality and $K\left(\mbf{X},\mbf{X}\right)$ is the so-called \emph{Gram} matrix that defines the covariance between the input values $\mbf{X}$ associated with $\mbf{y}$ through a kernel $k\left(\xvec_p,\xvec_q\right)$:
\begin{equation}
K_{pq} = k\left(\xvec_p,\xvec_q\right)=\noise_\mrm{f}\exp\left(-\frac{1}{2\ell^2}\norm{\xvec_p-\xvec_q}^2\right)
\label{EQrbfkernel}
\end{equation}
\noindent where instead of first defining a prior over weights and deriving an equivalent kernel as in \cref{EQlinmodelkernelform} we directly assume a kernel --- in this case the \emph{squared exponential} kernel \cite{Rasmussen} defined by a variance $\noise_\mrm{f}$ and a length scale $\ell$ --- which determines the level of correlation between points in input space. It can be shown that adopting the squared exponential kernel is equivalent to formulating a parametric model with an infinite number of basis functions $\bs\phi$ or a Bayesian neural network with one hidden layer with an infinite number of neurons \cite{Bishop,Rasmussen}.

With this definition for the covariance structure between function values, we can obtain the joint distribution of training targets $\mbf{t}_\obs$ by marginalizing (averaging) over all possible values of $\mbf{y}_\obs$:
\begin{equation}
\prob{\mbf{t}_\obs}=\int\prob{\mbf{t}_\obs\vert\yvec_\obs}\prob{\yvec_\obs}\mrm{d}\yvec=\mathcal{N}\left(\mbf{t}_\obs\vert\mbf{0},{K}\left(\mbf{X}_\obs,\mbf{X}_\obs\right)+\noise_\mrm{n}\mbf{I}\right)
\label{EQgptargetprior}
\end{equation}
\noindent which allows us to incorporate information from \dataset\ by defining a joint distribution between training values and new predictions for the function value:
\begin{equation}
\prob{\begin{bmatrix}\mbf{t}_\obs\\\mbf{y}_*\end{bmatrix}}=\mathcal{N}\left(\begin{bmatrix}\mbf{t}_\obs\\\mbf{y}_*\end{bmatrix}\Bigg\vert\mbf{0},\begin{bmatrix}K\left(\mbf{X}_\obs,\mbf{X}_\obs\right)+\noise_\mrm{n}\mbf{I}&K\left(\mbf{X}_\obs,\mbf{X}_*\right)\\K\left(\mbf{X}_*,\mbf{X}_\obs\right)&K\left(\mbf{X}_*,\mbf{X}_*\right)\end{bmatrix}\right) 
\label{EQgpjointprior}
\end{equation}
\noindent and subsequently find the conditional probability $\prob{\mbf{y}_*\vert\mbf{t}_\obs}$ of the new function values given that values coming from \dataset\ are known:
\begin{equation}
\prob{\mbf{y}_*\vert\mbf{t}_\obs} = \mathcal{N}\left(\mbf{y}_*\vert\mbf{m}_*,\mbf{S}_*\right)
\label{EQgpposterior}
\end{equation}
\noindent which is again a Gaussian distribution with mean and covariance given by \cite{Bishop,Rasmussen}:
\begin{equation}
\mbf{m}_* = K\left(\Xvec_*,\Xvec_\obs\right)\left(K\left(\Xvec_\obs,\Xvec_\obs\right)+\noise_\mrm{n}\mbf{I}\right)^{-1}\mbf{t}_\obs
\label{EQgpmstar}
\end{equation}
\begin{equation}
\mbf{S}_* = K\left(\Xvec_*,\Xvec_*\right)-K\left(\Xvec_*,\Xvec_\obs\right)\left(K\left(\Xvec_\obs,\Xvec_\obs\right)+\noise_\mrm{n}\mbf{I}\right)^{-1}K\left(\Xvec_\obs,\Xvec_*\right)
\label{EQgpsstar}
\end{equation}

Defining $\mbf{k}_*=K\left(\Xvec_\obs,\xvec_*\right) \in \real^{N\times 1}$ and $\mbf{K}_\obs=K\left(\Xvec_\obs,\Xvec_\obs\right) \in \real^{N\times N}$, we can arrive at shorter expressions for the expectation and variance of the function value at a single new point $\xvec_*$:
\begin{equation}
\mathbb{E}\left[y_*\vert\xvec_*\right] = \mbf{k}_*^\T\left(\mbf{K}_\obs+\noise_\mrm{n}\mbf{I}\right)^{-1}\mbf{t}_\obs
\label{EQgppredaverage}
\end{equation}
\begin{equation}
\mathbb{V}\left[y_*\vert\xvec_*\right] = k\left(\xvec_*,\xvec_*\right) - \mbf{k}_*^\T\left(\mbf{K}_\obs+\noise_\mrm{n}\mbf{I}\right)^{-1}\mbf{k}_*
\label{EQgppredvariance}
\end{equation}
\noindent Note that here we opt for directly using the function value $y_*$ instead of its noisy version $t_*$ to define our surrogate models. This effectively makes the adaptive components of the acceleration framework, which are driven by increases in the variance given by \cref{EQgppredvariance}, less sensitive to changes in $\sigma_\mrm{n}^2$. 

\begin{figure}
\centering
\begin{subfigure}[b]{0.45\textwidth}
\includegraphics[scale=0.9]{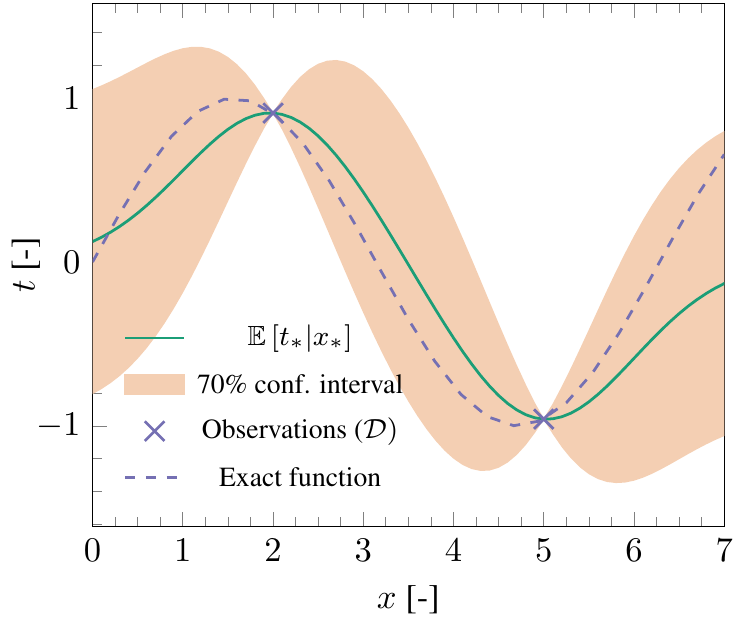}
\caption{Only target observations}
\end{subfigure}
\begin{subfigure}[b]{0.45\textwidth}
\includegraphics[scale=0.9]{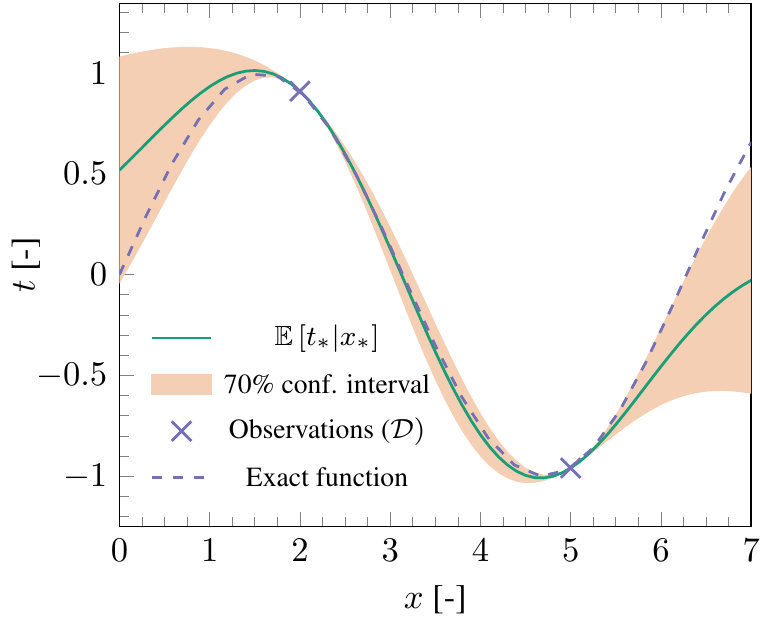}
\caption{Target and derivative observations}
\end{subfigure}
\caption{GP predictions for the function $t=\sin(x)$ constructed with two observations. Inclusion of derivative observations improves predictions around training points and reduces uncertainty.}
\label{FIgpdemo}
\end{figure}

In \cref{FIgpdemo}a we demonstrate a GP model by plotting predictions based on two observations of the noiseless target function $t=\sin(x)$. The conditioning of \cref{EQgpposterior} guarantees that $y_*$ coincides with the observations at the training points $\mbf{X}_\obs$. Away from the training space \dataset\ the GP returns to its zero-mean prior with high variance. This behavior is desirable for our application, since this uncertainty can be used as a trigger to either making new observations at the existing full-order integration points or adding new points to the full-order set when necessary. 

\subsection{Predicting derivatives and including derivative observations}

In building constitutive model surrogates, in addition to predicting new function values (\cref{EQgppredaverage}), their derivatives with respect to the input $\xvec_*$ are also needed. These can be computed by differentiating \cref{EQgppredaverage}: 
\begin{equation}
\frac{\partial}{\partial\xvec_*}\mathbb{E}\left[y_*\vert\xvec_*\right]=\sum_i^Na_i\frac{\partial}{\partial\xvec_*}k\left(\xvec_*,\xvec_i\right)
\label{EQderivs}
\end{equation}
\noindent where $a_i$ is the $i$-th component of the vector $\mbf{a}=\left(\mbf{K}+\noise_\mrm{n}\mbf{I}\right)^{-1}\mbf{t}$ and the derivative of the kernel is given by:
\begin{equation}
k'\left(\xvec_p,\xvec_q\right)\equiv\frac{\partial}{\partial\xvec_p}k\left(\xvec_p,\xvec_q\right) = -\frac{1}{\ell^2}\left(\xvec_p-\xvec_q\right)k\left(\xvec_p,\xvec_q\right)
\label{EQfstderivkernel}
\end{equation}

However, since we can also observe these derivatives (in the form of $\mbf{D}^\Omega$ from \cref{EQfe2stress}) at the anchor models, it is also interesting to include this information in the GP in order to improve its predictions. This is achieved by recognizing that the derivative of a GP remains Gaussian and defining the covariance between function values and derivatives as \cite{Solak2003}:
\begin{equation}
\mrm{cov}\left(\mbf{y}'_p,y_q\right)=\frac{\partial}{\partial\xvec_p}k\left(\xvec_p,\xvec_q\right)
\quad\quad
\mrm{cov}\left(\mbf{y}'_p,\mbf{y}'_q\right) = \frac{\partial^2}{\partial\xvec_q\partial\xvec_p}k\left(\xvec_p,\xvec_q\right) + \sigma^2_\mrm{d}
\label{EQderivcovariances}
\end{equation}
\noindent where $\sigma^2_\mrm{d}$ is a noise parameter that represents the uncertainty associated to observing the derivatives, the first-order derivative of the kernel is given in \cref{EQfstderivkernel} and its second derivative is a matrix given by:
\begin{equation}
k''\left(\xvec_p,\xvec_q\right)\equiv\frac{\partial^2}{\partial\xvec_q\xvec_p}k\left(\xvec_p,\xvec_q\right) = \frac{1}{\ell^2}\left(\mbf{I}-\frac{1}{\ell^2}\left(\xvec_p-\xvec_q\right)\left(\xvec_p-\xvec_q\right)^\T\right)k\left(\xvec_p,\xvec_q\right)
\label{EQsndderivkernel}
\end{equation}

With these definitions, the joint prior distribution of training targets becomes \cite{Pruher2016}:
\begin{equation}
\prob{\begin{bmatrix} \mbf{t}_\obs\\\mbf{t}'_\obs\end{bmatrix}}
= \mathcal{N}\left(\begin{bmatrix} \mbf{t}_\obs\\\mbf{t}'_\obs\end{bmatrix}\Bigg\vert\mbf{0},\begin{bmatrix}\mbf{K}_\obs+\noise_\mrm{n}\mbf{I} & \mbf{K}_\mrm{td}\\\mbf{K}_\mrm{td}^\T&\mbf{K}_\mrm{dd} \end{bmatrix}\right) \equiv \mathcal{N}\left(\begin{bmatrix} \mbf{t}_\obs\\\mbf{t}'_\obs\end{bmatrix}\Bigg\vert\mbf{0},\overline{\mbf{K}}_\obs\right)
\label{EQgptargetpriorderivs}
\end{equation}
\noindent where $\mbf{K}_\obs$ is the same kernel matrix appearing in \cref{EQgppredaverage} and $\mbf{K}_\mrm{td}\in\real^{N\times ND}$ and $\mbf{K}_\mrm{dd}\in\real^{ND\times ND}$ are composed of blocks of $k'\left(\xvec_p,\xvec_q\right)$ and $k''\left(\xvec_p,\xvec_q\right)$, respectively, with $\overline{\mbf{K}}_\obs$ being the resultant $N(D+1)\times N(D+1)$ covariance matrix ($D$ is the dimensionality of \xvec). Making point predictions is done in a similar way as in \cref{EQgppredaverage}:
\begin{equation}
\overline{\mbf{k}}_* = \begin{bmatrix} \mbf{k}_* & k'\left(\xvec_1,\xvec_*\right) & \cdots & k'\left(\xvec_N,\xvec_*\right)\end{bmatrix}^\T
\end{equation}
\begin{equation}
\mathbb{E}\left[y_*\vert\xvec_*\right] = \overline{\mbf{k}}_*^\T\overline{\mbf{K}}_\obs^{-1}\overline{\mbf{t}}_\obs\quad\quad\mathbb{V}\left[y_*\vert\xvec_*\right] = k\left(\xvec_*,\xvec_*\right) - \overline{\mbf{k}}_*^\T\overline{\mbf{K}}_\obs^{-1}\overline{\mbf{k}}_* + \noise_\mrm{n}
\end{equation}
\noindent where the target vector now includes the observed derivatives:
\begin{equation}
\overline{\mbf{t}}_\obs = \begin{bmatrix} t_{\obs 1}&\cdots&t_{\obs N}&\mbf{t}'_{\obs 1}&\cdots&\mbf{t}'_{\obs N}\end{bmatrix}^\T
\end{equation}
\noindent and predictions for the tangent are done as in \cref{EQderivs} but now by differentiating $\overline{\mbf{k}}_*$ instead of $\mbf{k}_*$. \cref{FIgpdemo} shows GP predictions for $t=\sin(x)$ with two training points with and without including derivative observations. Now the conditioning not only constrains the predictions to agree with the target values in \dataset\ but also with their derivatives. This makes effective use of the limited amount of information coming from a small number of \emph{online} observations, as making a single gradient observation is equivalent to adding $D$ extra function observations around the target value.

\subsection{Hyperparameter optimization}
\label{SEhyperparameteroptimization}

The process variance $\noise_\mrm{f}$ and length scale $\ell$ that compose the kernel and the target noise $\noise_\mrm{n}$ are hyperparameters that should be learned from the dataset \dataset.  Since a full Bayesian treatment for these parameters --- introducing a prior, deriving a posterior and marginalizing --- usually demands the use of expensive numerical techniques such as Markov Chain Monte Carlo (MCMC), here we opt for a maximum likelihood solution in order to minimize the computational overhead associated with the \emph{online} calibration of the GP models. Furthermore, due to the limited amount of data available for estimation, here we refrain from optimizing for the derivative noise $\sigma^2_\mrm{d}$ and instead assume derivative observations are noiseless.

The aim here is to maximize the marginal likelihood, obtained by averaging the probability that the model reproduces the training targets over all possible values of $\overline{\mbf{y}}_\obs$ associated with $\overline{\mbf{t}}_\obs$:
\begin{equation}
\prob{\overline{\mbf{t}}_\obs\vert\noise_\mrm{f},\noise_\mrm{n},\ell}=\int \prob{\overline{\mbf{t}}_\obs\vert\overline{\mbf{y}}_\obs,\noise_\mrm{n}}\prob{\overline{\mbf{y}}_\obs\vert\noise_\mrm{f},\ell}\mrm{d}\overline{\mbf{y}}_\obs
\label{EQgpmarginallikelihood}
\end{equation}
\noindent which is a function that depends only on the hyperparameters. We optimize this marginal likelihood with a BFGS algorithm \cite{Fletcher1987}, which requires the computation of the gradient of $\prob{\overline{\mbf{t}}}$. Taking the natural logarithm of both sides of \cref{EQgpmarginallikelihood} and differentiating with respect to $\overline{\mbf{K}}_\obs$ yields an expression for the gradient:
\begin{equation}
\nabla\ln\prob{\overline{\mbf{t}}_\obs}=-\frac{1}{2}\tr\left(\overline{\mbf{K}}^{-1}_\obs\nabla\overline{\mbf{K}}_\obs\right) + \frac{1}{2}\overline{\mbf{t}}_\obs^\T\overline{\mbf{K}}^{-1}_\obs\nabla\overline{\mbf{K}}_\obs\overline{\mbf{t}}_\obs
\end{equation}
\noindent and $\nabla\overline{\mbf{K}}_\obs$ can be obtained in a straightforward manner by differentiating \cref{EQrbfkernel,EQfstderivkernel,EQsndderivkernel} with respect to each hyperparameter and reassembling the matrix as in \cref{EQgptargetpriorderivs}.

\section{Surrogate modeling framework}
\label{SEsurrogatemodelingframework}

\RestyleAlgo{ruled}
\SetKw{Continue}{\textbf{continue}}
\SetSideCommentRight
\SetKwIF{If}{ElseIf}{Else}{if}{:}{else if}{else}{endif}
\SetKwFor{For}{for}{:}{}
\SetKwFor{While}{while}{:}{}
\SetKwInput{Params}{Params}

In this section, we use the techniques of \cref{SEbayesianparametricregression} to build an adaptive surrogate modeling framework for \fetwo. Following \cref{FIintro}, we introduce a surrogate model $\mathcal{S}$ trained on a set of constitutive observations \dataset\ obtained from a small number of fully-solved \emph{anchor} models (gathered in the set \Ful) subjected to the strain histories from their respective anchoring integration points. For bulk homogenization, we define $\mathcal{S}$ as:
\begin{equation}
\bs\sigma^\Omega\left(\bs\varepsilon^\Omega,\dataset\right) = \mbf{D}^\Omega_\mrm{e}\bs\varepsilon^\Omega + \mathbb{E}\left[\widehat{\bs\sigma}^\Omega\left(\bs\varepsilon^\Omega,\dataset\right)\right]
\label{EQcorrectionmodel}
\end{equation}
\noindent where $\mbf{D}^\Omega_\mrm{e}$ is an initial stiffness obtained from the first computed integration point and $\widehat{\bs\sigma}$ is a stress correction given by the mean GP response, with the variance being used exclusively for adaptivity purposes. The consistent tangent stiffness is the combination of $\mbf{D}^\Omega_\mrm{e}$ and the derivatives of the expected stress corrections (\cref{EQderivs}).

We opt for modeling each component of $\widehat{\bs\sigma}$ (three components for the 2D examples treated here, with which we implicitly enforce symmetry to the stress tensor) independently, each with its own zero-mean GP model. This strategy implicitly states that components of the stress tensor are a priori uncorrelated, which greatly simplifies the structure of the covariance matrix of the GP models and reduces the number of hyperparameters to be estimated. Furthermore, we assume for now that $\widehat{\bs\sigma}$ depends only on the current strain value $\bs\varepsilon^\Omega$, which means that the GP model loads, unloads and reloads along the same path. This is a limitation to be addressed in future versions of the framework that currently hinders its ability to treat general non-monotonic load paths\footnote{This can be addressed by augmenting the input space, the most straightforward way being using both $\bs\varepsilon^\Omega$ and $\Delta\bs\varepsilon^\Omega$ as inputs to the GP models.}. It is worth mentioning, however, that path dependency is already partially accounted for since different strain histories will lead to the construction of different surrogates. 

Note that the generality of the $\mathcal{S}$ model given by \cref{EQcorrectionmodel} is not compromised by the elastic-correction additive decomposition since $\widehat{\bs\sigma}$ can take any shape, but for initially linear-elastic materials this split improves the robustness of the surrogate in two ways. Firstly, it helps the GP in reproducing the initial elastic behavior at a moment when the surrogate models have very little data to work with. Secondly, it aids in preventing the occurrence of spurious strain localization as the GP moves back to its zero prior: away from the training points, the model of \cref{EQcorrectionmodel} returns instead to a linear-elastic response.

\begin{figure}
\centering
\includegraphics[width=0.7\textwidth]{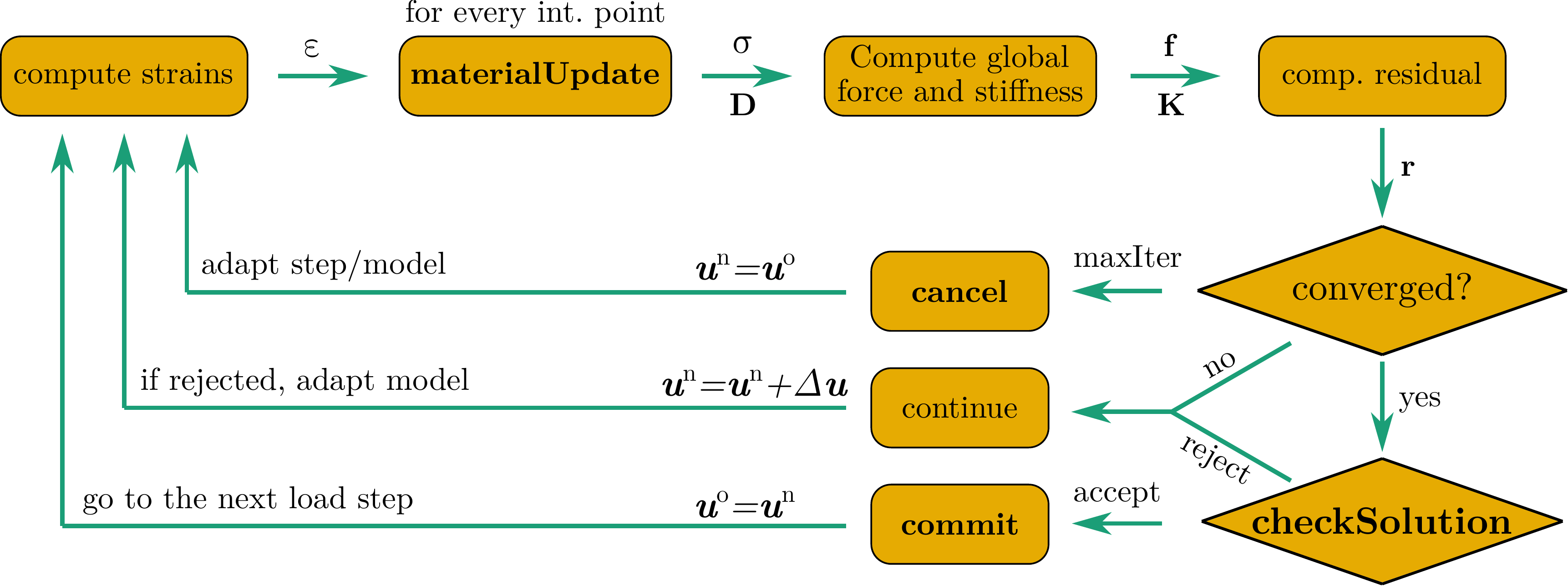}
\caption{Schematic analysis flow for finite element problems involving nonlinear material behavior. The subscripts $o$ and $n$ refer to old (converged) and new (current) values, respectively. The proposed acceleration framework focuses on the steps marked in bold.}
\label{FIflowchart}
\end{figure}

In order to position the present active learning approach within a general FEM implementation, it is useful to recall the main steps involved in finding equilibrium solutions for nonlinear finite element models. These are shown in \cref{FIflowchart}. A solution for $\mbf{u}$ is obtained iteratively by minimizing a global force residual $\mbf{r}$ computed from the material response at every integration point (\texttt{materialUpdate}). Upon convergence and before moving to the next time step, the current solution is checked (\texttt{checkSolution}) and can be rejected if there is a need to adapt the model (\eg nucleate/propagate cracks, refine the mesh, change constitutive models). Once the solution is converged and accepted, material history is updated (\texttt{commit}) and the model moves to the next time step. 

In an \fetwo\ model, \cref{FIflowchart} can represent the macroscopic solution loop and the micromodels can be seen as a material-like entity, which leads to another similar solution loop being embedded within the \texttt{materialUpdate} routine. Here we exploit this modularity and discuss the implementation of the learning approach as a surrogate to an arbitrary material model denoted as \texttt{fullModel}. We can therefore implement the approach as a material wrapper that encapsulates \texttt{fullModel} and handles learning and prediction tasks. \cref{ALupdate,ALcheckcommit,ALcommit,ALcancel} show how the material routines marked in bold in \cref{FIflowchart} are implemented for this wrapper. In the following, we elaborate on each of these components and present the main implementational aspects of the framework in detail.

\subsection{Initial sampling}
\label{SEinitialsampling}

Since there is no \emph{offline} training, the surrogate model starts with no prior information on the constitutive behavior being approximated ($\dataset=\varnothing$) and no \emph{anchor} models initially present ($\Ful=\varnothing$). We therefore introduce an initialization step by solving the first time increment under the assumption that all integration points are deforming elastically in order to obtain an initial strain distribution which we use to initialize \Ful\ and \dataset. During this step, we use the fully-solved model only once in order to obtain the initial stiffness $\mbf{D}_\mrm{e}$ and we assume all other points have the same stiffness in order to avoid further full-order computations (\cref{ALupdate}). As we will see shortly, this is not a limiting assumption since this first elastic approximation is rejected and the first time step is revisited after the GP models and initial anchors are initialized.

\begin{algorithm}
\caption{The \texttt{materialUpdate} routine}
\label{ALupdate}
\BlankLine
\KwIn{strain $\strain$ at the integration point $p$}
\KwOut{stress $\stress$ and stiffness $\mbf{D}$ at the integration point}
\BlankLine
\uIf{initialization step}{
\uIf{$\mbf{D}_\mrm{e}$ is not initialized (\ie $p$ is the first point ever computed)}{
use original constitutive model:$\left(\stress,\mbf{D}\right)$ $\leftarrow$ $\mathtt{fullModel::materialUpdate}\left(\strain\right)$\;
initialize elastic stiffness: $\mbf{D}_\mrm{e}\leftarrow\mbf{D}$\;
}
\Else{
assume $p$ has the same stiffness as the first computed point: $\bs\sigma\leftarrow\mbf{D}_\mrm{e}\bs\varepsilon;\,\,\mbf{D}\leftarrow\mbf{D}_\mrm{e}$\;
}
}
\Else{
compute elastic stresses: $\stress_\mrm{e}\leftarrow\mbf{D}_\mrm{e}\strain$\;
use GP to predict a constitutive correction: 
$\left(\mathbb{E}\left[\widehat{\stress},\widehat{\mbf{D}}\right],\mathbb{V}\left[\widehat{\stress}\right]\right)\leftarrow\mathtt{GP}\left(\strain\right)$\;
approximate the response: $\stress\leftarrow\stress_\mrm{e}+\mathbb{E}\left[\widehat{\stress}\right],\quad\mbf{D}\leftarrow\mbf{D}_\mrm{e}+\mathbb{E}\left[\widehat{\mbf{D}}\right]$\;
compute the uncertainty indicator $\gamma$ as in \cref{EQgamma}\;
\If{$\gamma>\gamma_\mrm{cancel}$}{\textbf{cancel time step (\cref{ALcancel})}\;}
\If{softening is detected ($D_{ii}<0$)}{penalize point uncertainty: $\gamma\leftarrow\gamma-D_{ii}$\;}
\If{last step has been cancelled}{switch to a modified Newton-Raphson strategy: $\mbf{D}\leftarrow\mbf{D}_\mrm{e}$\;}
store values for this point: $\strain^\mrm{n}_p\leftarrow\strain;\,\,\gamma_p^\mrm{n}\leftarrow\gamma$\;
}
\Return{$\bs\sigma,\mbf{D}$}
\end{algorithm}

After this first approximation for the solution is obtained, the \texttt{checkSolution} routine of \cref{ALcheckcommit} is called. Since the strains at every integration point are stored, an informed initial choice for \Ful\ can be made by clustering the integration points using a $k$-means clustering algorithm \cite{Lloyd1982}. For each cluster, the point closest to the cluster centroid is chosen, added to \Ful\ and immediately sampled. The GP models are then prepared to make predictions by computing and storing factorized versions of their covariance matrices. 

\begin{algorithm}
\caption{The \texttt{checkSolution} routine}
\label{ALcheckcommit}
\BlankLine
\uIf{first time step}{
initialize dataset and anchor point set : $\dataset\leftarrow\varnothing$, $\Ful\leftarrow\varnothing$\;
divide the integration points into $k$ clusters in strain space: $\left(\bs{\mathcal{C}},\overline{\strain}^\mrm{n}\right)\leftarrow\mathtt{KMC}\left(\strain^\mrm{n},k\right)$\;
\For{every cluster $\mathcal{C}_i\in\bs{\mathcal{C}}$}{
find representative point $p=\underset{j\in\mathcal{C}_i}{\mathrm{arg\,min}}\norm{\strain^\mrm{n}_{j}-\overline{\strain}^\mrm{n}_i}$\;
anchor a fully-solved model at point $p$: $\Ful\leftarrow\Ful\cup p$\;
solve the newly-created full model once: $\left(\stress,\mbf{D}\right)$ $\leftarrow$ $\mathtt{fullModel::materialUpdate}\left(\strain^\mrm{n}_p\right)$\;
compute corrections and add data: $\widehat{\stress}\leftarrow\stress-\mbf{D}_\mrm{e}\strain^\mrm{n}_p;\,\,\widehat{\mbf{D}}\leftarrow\mbf{D}-\mbf{D}_\mrm{e};\,\,\dataset\leftarrow\dataset\cup\left(\strain^\mrm{n}_{{p}},\left[\widehat{\stress}\quad\widehat{\mbf{D}}\right]\right)$\;
}
update GP and store marginal likelihood: $\mathtt{updateGP}\left(\mathcal{D}\right);\,\,L\leftarrow\ln p\left(\mbf{y}\vert\mbf{x}\right)$\;
\textbf{reject solution (revisit the first time step)}\;
}
\Else{
update $\mathcal{D}$ based on the tolerance $\gamma_\mrm{tol}$ using \cref{ALadddata}: $\mathcal{D}\leftarrow\mathcal{D}\cup\mathtt{addData}\left(\gamma_\mrm{tol}\right)$\;
\uIf{\added\ has not changed}{
\textbf{accept solution as is (\cref{ALcommit})}\;}
\Else{
refactor covariance matrix: $\mathtt{updateGP}\left(\mathcal{D}\right)$\;
\If{$\norm{L / \ln p\left(\overline{\mbf{t}}\right)} > L_\mrm{retrain}$}{recompute GP hyperparameters: $\mathtt{retrainGP}\left(\mathcal{D}\right);\,\,L\leftarrow\ln p\left(\overline{\mbf{t}}\right)$\;}
\textbf{reject solution (continue on the same time step)}\;
}
}
\end{algorithm}

Initial values for the hyperparameters $\noise_\mrm{f}$, $\noise_\mrm{n}$ and $\ell$ can either be provided by the user based on prior knowledge (\eg from using the same material on other models)\footnote{Hyperparameters reflect intrinsic patterns in the constitutive manifold being inferred and should therefore be fairly insensitive to the macroscopic model being solved.} or be estimated \emph{online} with an initial training procedure. We do this here by creating a fictitious anchor model for each cluster and loading it monotonically in the strain direction seen by the point closest to the cluster centroid. In order to avoid redundancy (\eg adding linear-elastic data coming from multiple models), we first initialize the GP with the standard hyperparameter values \sigf{1.0}, \length{1.0e-2}, \sign{0.0} and use a variance threshold as a dosing mechanism to ensure the added data is uniformly spaced. After optimizing for the hyperparameters, we discard the data obtained from the fictitious anchors since there is no guarantee any of the real integration points will follow the same strain history. Once \dataset\ and \Ful\ are initialized, the first time step is then revisited by rejecting the initial elastic approximation done in \cref{ALupdate}. Note that this happens only once at the beginning of the analysis.

\subsection{Model adaptivity}
\label{SEmodeladaptivity}

After the initialization step, the GP surrogates are used to compute $\widehat{\bs{\sigma}}$ and the response is approximated as in \cref{EQcorrectionmodel}. Since an independent GP model is used for each stress component, we adopt a single uncertainty indicator $\gamma$ given by:
\begin{equation}
\gamma = \max_i^{n_\sigma}\left(\sqrt{\mathbb{V}\left[\widehat{\sigma}_i\right]}\right)
\label{EQgamma}
\end{equation}
\noindent with $n_\sigma$ being the number of stress components, and use it to drive the adaptive components of the framework\footnote{The indicator $\gamma$ used here is expressed in units of stress and is therefore an absolute measure. We also experimented with a number of relative indicators but found those to be excessively sensitive to changes in the hyperparameter values.}. 

Values of $\gamma$ for all integration points are updated at every global Newton-Raphson iteration (\cref{ALupdate}) and their final values $\bs{\gamma}^\mrm{n}$ (upon global convergence) are used to drive a greedy refinement of the GP models (\cref{ALcheckcommit}) and ensure the surrogate model remains accurate. This is enforced by the tolerance parameter $\gamma_\mrm{tol}$ through the \texttt{addData} routine of \cref{ALadddata}. Here we define an additional set \added\ of tracked anchors that contains the models for which data has been added during the present time step. The set \added\ is used both to avoid repeatedly adding data from a single anchor model at any given time step and to make sure the single added data point is updated to reflect the latest equilibrium solution (since from \cref{FIflowchart}, multiple $\mathtt{materialUpdate}\rightarrow\mathtt{checkSolution}\rightarrow\mathtt{continue}$ cycles can occur before history is committed).

The uncertainty level $\gamma$ is checked at every integration point and the anchor model associated with the point having the highest value is chosen for sampling. In order to keep the number of models in \Ful\ to a minimum, \texttt{addData} gives priority to models already in \Ful\ and only considers other potential anchoring points if every model in \Ful\ with an uncertainty higher than $\gamma_\mrm{tol}$ has already been sampled on the current time step. Since the macroscopic problem being solved imposes a degree of similarity between integration points, the idea is that adding data from points in \Ful\ should also help reduce the uncertainty of nearby integration points. This can therefore be seen as a \emph{greedy} data selection approach \cite{Rasmussen} that aims at keeping both \Ful\ and \dataset\ as small as possible. If all models in \Ful\ have already been sampled and a new anchoring point must be chosen, we first recover the material history of the newly-created model by making it revisit the complete cumulative strain history $\strain_\mrm{h}^{\overline{p}}$ of its anchoring point $\overline{p}$. Also note that we avoid adding data from models in \Ful\ undergoing unloading or reloading. This is a consequence of assuming a unique relationship between strains and stresses: data coming from anchor models unloading through a different path would be erroneously interpreted by the GP models as noise.

\begin{algorithm}
\caption{The \texttt{addData} routine}
\label{ALadddata}
\BlankLine
\KwIn{A tolerance $\gamma_\mrm{tol}$}
\KwOut{A single data point to be included in $\mathcal{D}$}
\BlankLine
update data sampled from \added\ to reflect the new equilibrium strains $\strain^\mrm{n}_{p\in\added}$\;
update GP models and uncertainty indicators: $\mathtt{updateGP}\left(\mathcal{D}\right);\,\,\bs{\gamma}^\mrm{n}\leftarrow\bs{\gamma}^{\strain^\mrm{n}}$\;
try to find an untracked anchor point to sample: $\overline{p}=\underset{p\in\Ful,p\notin\added}{\mrm{arg\,max}}\,\,\gamma^\mrm{n}_p\vert\gamma^\mrm{n}_p>\gamma_\mrm{tol}$\;
\If{$\overline{p}$ cannot be found}{
try to find a new anchoring point: $\overline{p}=\underset{p\notin\Ful,p\notin\added}{\mrm{arg\,max}}\,\,\gamma^\mrm{n}_p\vert\gamma^\mrm{n}_p>\gamma_\mrm{tol}$\;
}
update list of tracked points: $\added\leftarrow\added\cup\overline{p}$\;
update list of anchor points: $\Ful\leftarrow\Ful\cup\overline{p}\,\,$ if $\,\,\overline{p}\notin\Ful$\;
call $\mathtt{fullModel::materialUpdate(\strain^\mrm{h}_{\overline{p}})}$ to gradually bring the anchor at $\overline{p}$ to the current time step\; 
compute model response at $\strain^\mrm{n}_{\overline{p}}$: $\left(\stress,\mbf{D}\right)$ $\leftarrow$ $\mathtt{fullModel::materialUpdate}\left(\strain^\mrm{n}_{\overline{p}}\right)$\;
compute corrections: $\widehat{\stress}^\mrm{n}_{\overline{p}}\leftarrow\stress-\mbf{D}_\mrm{e}\strain^\mrm{n}_{\overline{p}};\,\,\widehat{\mbf{D}}^\mrm{n}_{\overline{p}}\leftarrow\mbf{D}-\mbf{D}_\mrm{e}$\;
\Return $\left(\strain^t_{\overline{p}},\begin{bmatrix} \widehat{\stress}^t_{\overline{p}} & \widehat{\mbf{D}}^t_{\overline{p}}\end{bmatrix}\right)$ with $t$ being the latest step on which $\overline{p}$ is not unloading/reloading\;
\end{algorithm}

Upon making changes to \dataset, the covariance matrices of all GP models are updated and refactored. The addition of a new observation also leads to a change in marginal likelihood (\cref{EQgpmarginallikelihood}). We therefore allow for a re-estimation of the hyperparameters to take place once the likelihood reaches a value $L_\mrm{retrain}$ times lower than the one computed after the latest estimation. If \added\ has not changed during the latest call to \texttt{addData}, the current solution is accepted as it is. Otherwise we continue with the same time step by rejecting the current solution, which will cause the global Newton-Raphson solver to keep searching for an equilibrium solution but now with updated GP models.


Two additional adaptivity safeguards are put in place in the \texttt{materialUpdate} routine of \cref{ALupdate}. Firstly, the time step can be canceled if the uncertainty is higher than a threshold $\gamma_\mrm{cancel}$ in order to avoid considering equilibrium solutions that are excessively far from the constitutive behavior being approximated. As we will discuss in \cref{SEsolutionrobustness}, canceling the current load step does not mean giving up on the analysis, with new solution attempts being made after including extra anchor models in \Ful\ and retraining the GP models. Secondly, the diagonal of the tangent stiffness matrix is checked for possible negative values, providing an indication of softening. This would indicate a switch from stresses to tractions is in order, but the decision for such a constitutive model switch should always be based on accurate information obtained from models in \Ful. Therefore, we flag the point for sampling by penalizing its uncertainty\footnote{In the examples of this paper we do not treat models for which such a switch from bulk to cohesive behavior would be necessary. The safeguard is therefore only triggered in rare occasions when the GP is trying to approximate a perfectly-plastic response.}. 

Apart from the very first update computed in order to obtain $\mbf{D}_\mrm{e}$, the expensive full-order model is never computed during \texttt{materialUpdate} (\cref{ALupdate}). An alternative to this approach would be to actually call the models in \Ful\ every time stresses at the anchoring points are computed. We do not opt for this alternative for two reasons. Firstly, using a single constitutive model (the surrogate $\mathcal{S}$) for the whole mesh avoids potential non-uniqueness issues that would arise, for instance, if points in \Ful\ are switching between different constitutive regimes (loading/unloading/softening). Secondly, refraining from constantly updating the models in \Ful\ results in significant gains in terms of acceleration by making the reduced model insensitive to the number of global Newton-Raphson iterations needed for convergence and allowing full-order models in \Ful\ to become dormant and essentially be removed from the analysis for as long as the uncertainty at the associated anchoring point does not increase. However, it is worth mentioning in passing that the more expensive alternative has the merit of allowing for an extra \emph{novelty detection} safeguard to be employed through which the deviation between predictions for $\stress$ coming from the full-order and surrogate models can be kept in check.

\subsection{Solution robustness}
\label{SEsolutionrobustness}

\cref{ALcheckcommit,ALadddata} ensure that only a single data point is added to \dataset\ at a time. This is done in order to avoid large perturbations to the current equilibrium solution caused by changes in the constitutive response\footnote{Adding data to the GP affects the response of all points within a hypersphere in strain space with a radius that depends on $\ell$. The resulting change in global response may cause the solver to diverge.}. Once a new observation is added, the current (now unconverged) solution is kept but the algorithm continues on the same time step until equilibrium is reestablished, at which point \texttt{checkSolution} is once again called. This process is repeated until $\gamma^\mrm{n}$ is lower than $\gamma_\mrm{tol}$ everywhere.  The resulting cycle of carefully adding data without significantly drifting away from equilibrium helps keeping the global solution scheme robust.

\begin{algorithm}
\caption{The \texttt{commit} routine}
\label{ALcommit}
\BlankLine
clear list of tracked models: $\added\leftarrow\varnothing$\;
store converged values: $\bs\gamma^\mrm{o}\leftarrow\bs\gamma^\mrm{n}$\;
update persistent strain history: $\strain^\mrm{h}\leftarrow\begin{bmatrix}\strain^\mrm{h} & \strain^\mrm{n}\end{bmatrix}$\;
\textbf{go to next time step}\;
\end{algorithm}

When the solution converges and $\gamma^\mrm{n}<\gamma_\mrm{tol}$ for all integration points, the \texttt{commit} routine is called (\cref{ALcommit}). Before moving to the next step, converged values for $\bs\gamma$ are updated. If the global solver fails to converge or if the additional uncertainty threshold $\gamma_\mrm{cancel}$ of \cref{ALupdate} is violated, the solution for the current time step is canceled (\cref{ALcancel}). Before making a new solution attempt, converged values are recovered and \cref{ALadddata} is used to sample the model with highest uncertainty even if its value for $\gamma$ is lower than $\gamma_\mrm{tol}$. Furthermore, we switch to a secant solution strategy by fixing the stiffness matrix to $\mbf{D}_\mrm{e}$ (\cref{ALupdate}). We empirically notice that, for certain combinations of hyperparameters, the surrogate model causes the global Newton-Raphson solver to lose robustness due to rapid changes in stiffness. Switching to a secant strategy after a canceled step assures a solution is found under this scenario, albeit with a lower convergence rate.

\begin{algorithm}
\caption{The \texttt{cancel} routine}
\label{ALcancel}
\BlankLine
recover converged values: $\bs\gamma^\mrm{n}\leftarrow\bs\gamma^\mrm{o}$\;
add data from the point with the highest variance (\cref{ALadddata}): $\mathcal{D}\leftarrow\mathcal{D}\cup\mathtt{addData}\left(0\right)$\;
refactor covariance matrix: $\mathtt{updateGP}\left(\mathcal{D}\right)$\;
\If{$\norm{L / \ln p\left(\overline{\mbf{t}}\right)} > L_\mrm{retrain}$}{recompute GP hyperparameters: $\mathtt{retrainGP}\left(\mathcal{D}\right);\,\,L\leftarrow\ln p\left(\overline{\mbf{t}}\right)$\;}
\textbf{go back to the beginning of the time step}\;
\end{algorithm}

\section{Numerical examples}
\label{SEexamples}

In this section, we put the proposed framework to the test on a number of numerical examples. The algorithms of \cref{SEsurrogatemodelingframework} have been implemented in an in-house Finite Element code using the open-source Jem/Jive C++ numerical analysis library \cite{jive}. We begin by demonstrating the applicability of the framework for \fetwo\ analysis and move on to performing an in-depth investigation of its performance with fast single-scale homogeneous examples that allow for extensive parametric studies to be performed. It is emphasized that, although it is our vision to use the surrogate model in a multiscale context, it does not matter for the purpose of testing the active learning framework whether the material update of the full model involves solving a micromechanical BVP or just evaluating a nonlinear constitutive relation.

\subsection{\fetwo\ demonstration}
\label{SEfe2demonstration}

\begin{figure}
\centering
\includegraphics[width=0.75\textwidth]{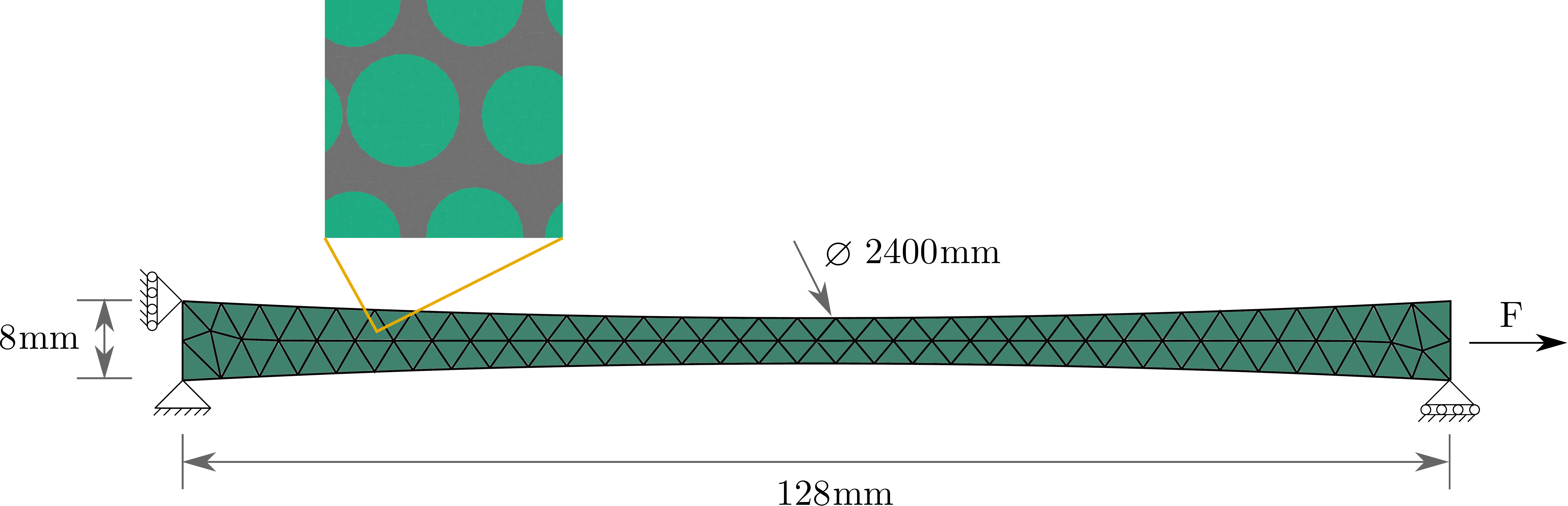}
\caption{\fetwo\ demonstration: A composite tapered bar loaded in transverse tension.}
\label{FIfe2dogmesh}
\end{figure}

The first example concerns a fiber-reinforced composite tapered specimen loaded in transverse tension. The geometry, mesh and boundary conditions are shown in \cref{FIfe2dogmesh}. A 4-fiber RVE model is embedded at each macroscopic integration point. The geometry and mesh of the micromodel are also shown in \cref{FIfe2dogmesh}. The fibers are modeled as linear elastic with properties $E=\SI{74000}{\mega\pascal}$ and $\nu=\SI{0.2}{}$. For the matrix we employ the pressure-dependent elastoplastic model proposed by Melro \etal \cite{Melro2013}, with $E=\SI{3130}{\mega\pascal}$, $\nu=\SI{0.37}{}$, $\nu_\mrm{p}=\SI{0.32}{}$ (plastic Poisson's ratio) and yield stresses given by:
\begin{equation}
\sigma_\tens = 64.80-33.6\e^{-\varepsilon^\p_\eq/0.003407}-10.21\e^{-\varepsilon^\p_\eq/0.06493}
\label{EQyieldsurface}
\end{equation}
\begin{equation}
\sigma_\comp = 81.00-42.0\e^{-\varepsilon^\p_\eq/0.003407}-12.77\e^{-\varepsilon^\p_\eq/0.06493}
\end{equation}
\noindent where $\varepsilon^\p_\eq$ is the equivalent plastic strain. The model is solved for 100 load steps, at which point the global macroscopic response is almost perfectly plastic and the strain localizes around the center of the specimen.

We run the reduced model with $k=1$, $\gamma_\mrm{tol}=\SI{0.3}{\mega\pascal}$, $\gamma_\mrm{cancel}=\SI{20}{\mega\pascal}$ and with hyperparameters estimated using a single micromodel loaded in the direction of the single $k$-means clustering centroid (in this case the average strain in the specimen). The analysis starts with $\Ful=\varnothing$ and $\dataset=\varnothing$ and ends with a total of 14 anchor models in \Ful\ (out of a total of 134 points) and $\vert\dataset\vert=73$ data points. 

\begin{figure}
\centering
\includegraphics[scale=1]{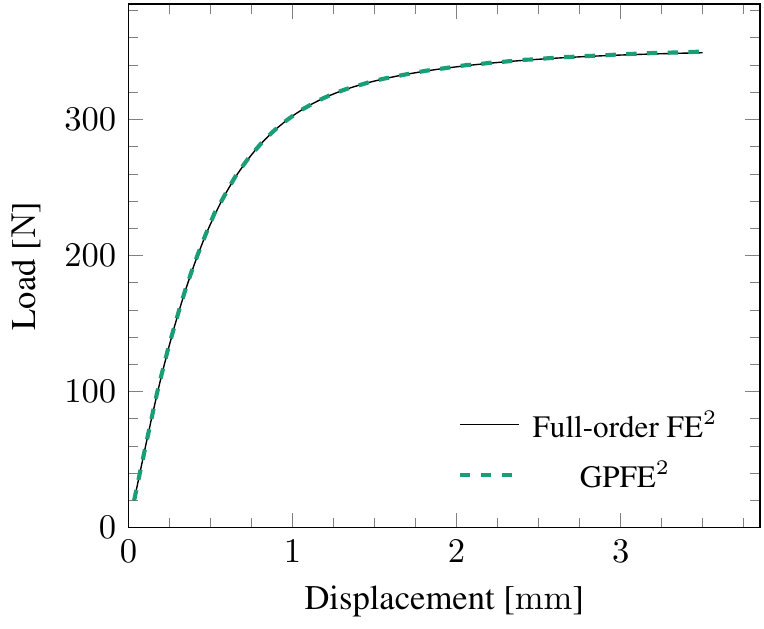}
\caption{Load-displacement curves obtained with the full-order \fetwo\ approach and with our reduction framework. The reduced model runs $\SI{22}{}$ times faster with only negligible loss of accuracy.}
\label{FIfe2doglodi}
\end{figure}

\begin{figure}
\centering
\includegraphics[scale=0.4]{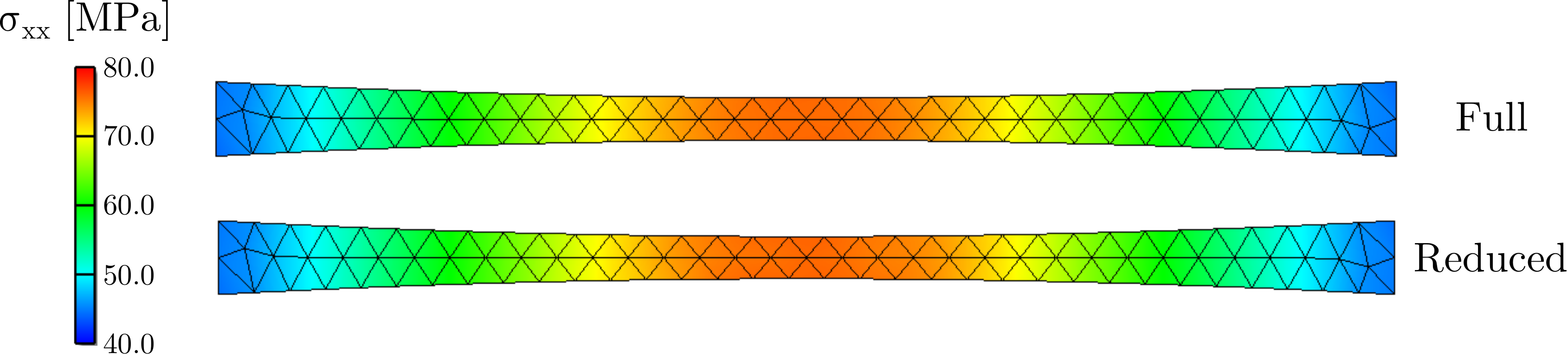}
\caption{\fetwo\ tapered bar: Stress fields at the end of the analysis for both full-order and reduced models. Our active learning approach predicts a stress distribution indistinguishable from the full-order one.}
\label{FIfe2dogsigxx}
\end{figure}

We plot the resultant load-displacement response for both full- and reduced-order models in \cref{FIfe2doglodi}. The active learning approach is able to capture the correct global model response with negligible loss of accuracy while running $22$ times faster than the full-order model. Of the \SI{35}{\second} execution time of the reduced model, a total of \SI{15}{\second} is spent on updating the GP models and using them for new predictions. For larger micromodels with denser meshes, the execution time of the reduced model will be dominated by computing the few micromodels included in \Ful\ and the overhead associated with the GP models will become negligible. For the full-order model, \SI{99}{\percent} of the execution time is spent solving the embedded micromodels, confirming the bottleneck assumption of \cref{EQfe2cost}. In \cref{FIfe2dogsigxx} we plot the horizontal stress distribution along the specimen at the last time step for both models. No discernible differences between the two stress distributions can be seen. 

\subsection{Performance and parametric sensitivity}

The implementation of \cref{SEsurrogatemodelingframework} allows for the active learning framework to supplant any full-order constitutive model. In the example of \cref{SEfe2demonstration}, this full-order model is the embedded RVE model of \cref{FIfe2dogmesh}. We now switch to a single-scale model with the homogeneous and elastoplastic material used in the previous example. This allows for an in-depth investigation on the performance of the reduction framework to be performed without loss of generality and without resorting to running a large number of expensive \fetwo\ simulations.

\begin{figure}
\centering
\includegraphics[width=0.7\textwidth]{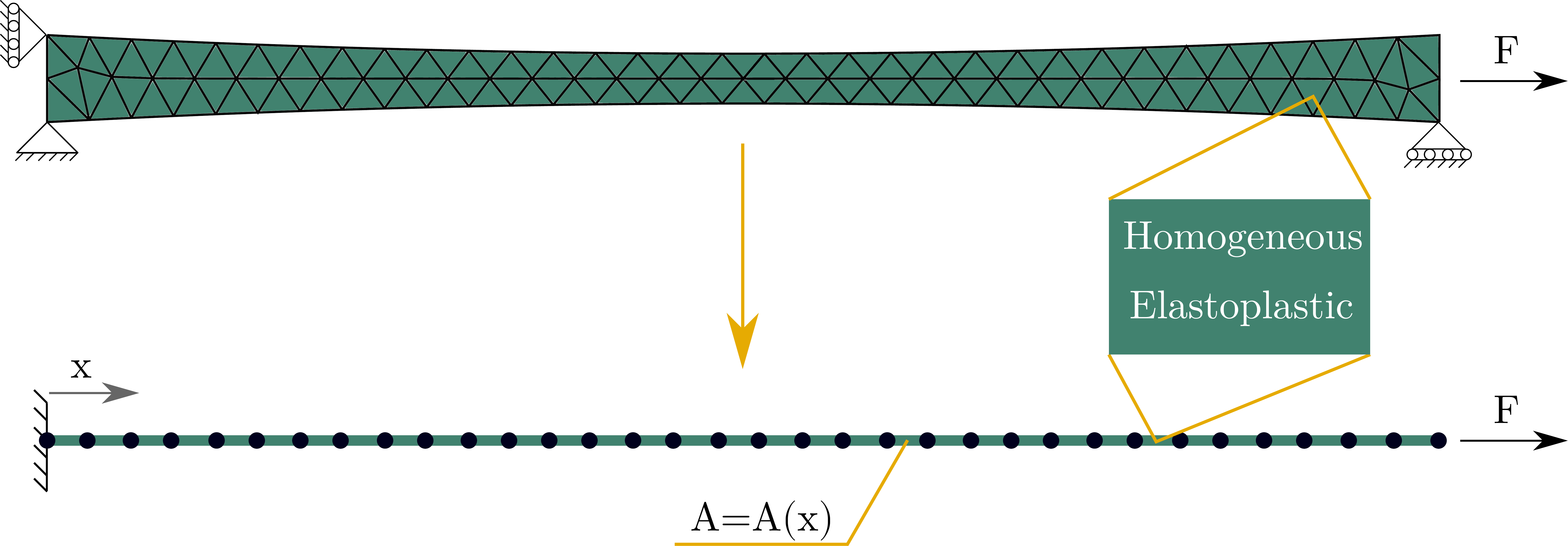}
\caption{Homogeneous models used to study the parametric sensitivity of the framework. The additional simplification to a bar problem allows for the constitutive space to be visualized in two dimensions ($\varepsilon_\mrm{xx}$-$\sigma_\mrm{xx}$).}
\label{FIbarmesh}
\end{figure}

We start the investigation by further simplifying the macroscopic model to the one-dimensional bar with variable cross-section area shown in \cref{FIbarmesh}, reducing the dimensionality of the constitutive space being approximated to only two dimensions ($\varepsilon_\mrm{xx}-\sigma_\mrm{xx}$). This change allows for the full model being inferred to be easily visualized and reduces the number of hyperparameters from nine for the two-dimensional case to only three\footnote{For each stress component we have three hyperparameters to be determined, namely $\noise_\mrm{f}$, $\noise_\mrm{n}$ and $\ell$.}. The original tapered geometry is simulated by making the cross-sectional area of the bar depend on the $x$ coordinate:
\begin{equation}
A(x) = 0.8 - 2.0\left(0.0534x - 0.000418x^2\right)
\end{equation}

\subsubsection{Evolution of the GP regression}

\begin{figure}
\centering
\includegraphics[width=0.75\textwidth]{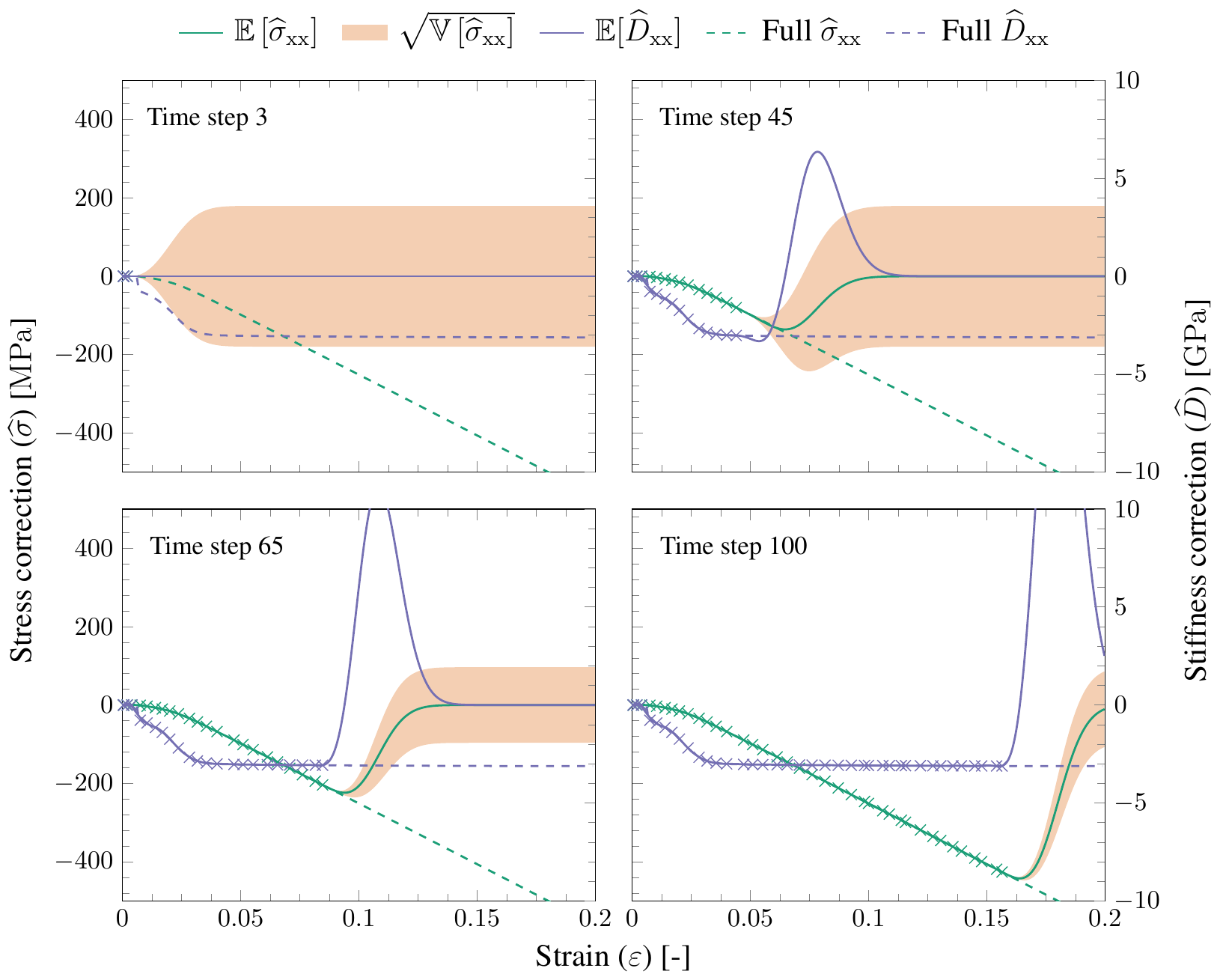}
\caption{Evolution of the GP constitutive model as data is gradually added. Dashed lines show the full-order response being approximated. The observations in \dataset\ are plotted as crosses.}
\label{FIbarevolution}
\end{figure}

We run the reduced model with $k=1$, \gtol{0.4}, $\gamma_\mrm{cancel}=\SI{20}{\mega\pascal}$ and $L_\mrm{retrain}=10$, with hyperparameter values obtained by loading a single material point in the $k$-means centroid direction. We can visualize the gradual improvement of the surrogate model as data is added by plotting the GP predictions together with the exact constitutive response being approximated at different moments throughout the analysis. This can be seen in \cref{FIbarevolution}, which shows snapshots made at four different time steps. 

As we start with an empty dataset, at first the GP model predicts linear-elastic behavior for the complete strain range, but with an uncertainty that quickly increases away from the observations (marked as crosses in \cref{FIbarevolution}). This increase is the mechanism that triggers the sampling of extra information from models in \Ful.  As more data is added to \dataset, the GP model is gradually refined and is able to reproduce the target full-order response with excellent accuracy, although we observe that even though the perfectly plastic response has a simple and constant shape, the GP is never able to extrapolate it. Away from the sampling points the GP predictions will always return to the prior with zero mean. By the end of the analysis, four models are present in \Ful\ and $\vert\dataset\vert=39$.

\subsubsection{Reduction efficiency}

\begin{figure}
\centering
\includegraphics[scale=1.0]{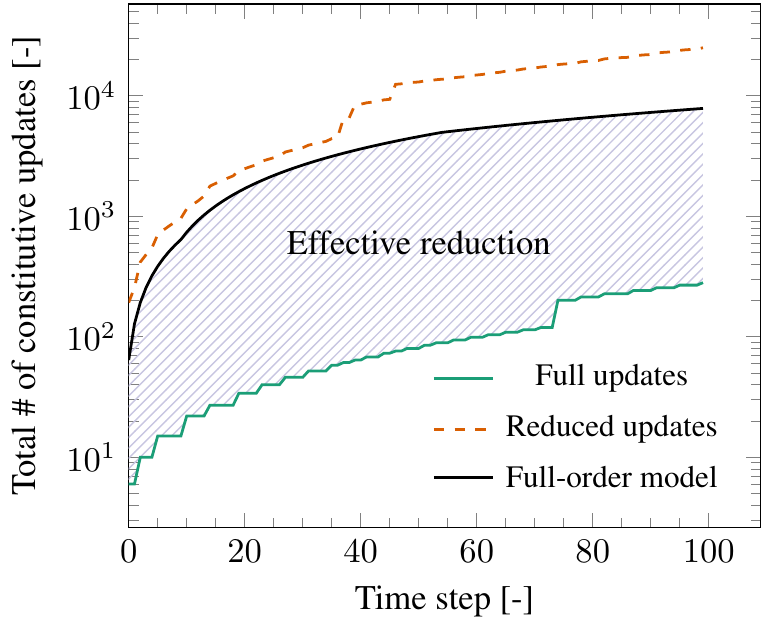}
\caption{Cumulative number of constitutive updates obtained with the full-order and reduced models. The shaded region represents the effective acceleration associated with the reduction framework.}
\label{FIbarrfactor}
\end{figure}

In \cref{SEconcurrentmultiscaleanalysis} we argued that the computational bottleneck of \fetwo\ lies on macroscopic constitutive model evaluations (\cref{EQfe2cost}). We can therefore have an indication of the acceleration promoted by the active learning approach simply by counting the number of full-order model evaluations. We plot in \cref{FIbarrfactor} the evolution of the cumulative number of constitutive updates for the example of the previous section together with the number of updates obtained by using the full-order model at every integration point. 

We see that the total number of material updates performed by the reduced model is significantly higher than that of the reference full-order model. Two different reasons for this increase can be identified. Firstly, four time step cancels occur after which the model switches to a secant approach that requires more iterations for convergence. Secondly, extra iterations are triggered every time a new observation is added to \dataset\ and the solution deviates from equilibrium (\cref{ALcheckcommit,ALadddata}).

However, given that the bottleneck assumption of \cref{EQfe2cost} holds, the effective acceleration brought by the adaptive reduction approach is only related to the number of times the anchor models are computed, and is therefore related to the shaded area shown in \cref{FIbarrfactor}. Here we define $R$ as the ratio between full model evaluations of the reference (full-order) and reduced-order models:
\begin{equation}
R=\frac{n_\mrm{ref}^\ful}{n_\mrm{GP}^\ful}
\label{EQreductionratio}
\end{equation}
\noindent and use it as a measure of acceleration. For the model of \cref{FIbarrfactor}, the ratio at the end of 100 time steps is $R=\SI{27.9}{}$.

\subsubsection{Influence of $\gamma_\mrm{tol}$}

The uncertainty tolerance $\gamma_\mrm{tol}$ is the main parameter controlling the active learning procedure: lower values of $\gamma_\mrm{tol}$ should lead to a higher sampling frequency while higher values should lead to smaller cardinalities for \dataset\ and \Ful\ at the cost of solution accuracy. In this section we put these claims to the test and investigate how much control can actually be exerted over the solution algorithm by changing $\gamma_\mrm{tol}$.

\begin{figure}
\centering
\includegraphics[scale=1.0]{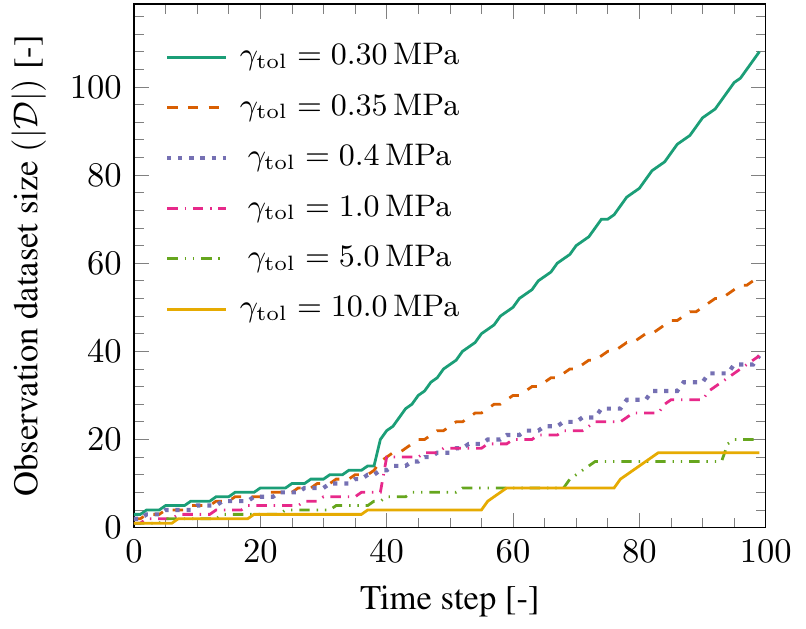}
\caption{Evolution of the cardinality of \dataset\ for different values of $\gamma_\mrm{tol}$. A stricter uncertainty tolerance leads to a high sampling frequency.}
\label{FIbardsize}
\end{figure}

Going back to the one-dimensional example of the previous section, we solve the problem for multiple values of $\gamma_\mrm{tol}$ between $\SI{0.3}{\mega\pascal}$ and $\SI{10.0}{\mega\pascal}$. In \cref{FIbardsize} we plot the evolution of the cardinality of the dataset \dataset\ for six different $\gamma_\mrm{tol}$ values. We see that $\gamma_\mrm{tol}$ indeed influences the number of observations added to \dataset, but only to a limited extent. This is due to the presence of the remaining adaptive components of the framework, namely the time step cancelling mechanism related to $\gamma_\mrm{cancel}$, the sampling of data with zero threshold after a canceled step and the hyperparameter retraining procedure linked to $L_\mrm{retrain}$. This combination of safeguards leads to datasets of similar sizes for most of the $\gamma_\mrm{tol}$ values adopted here. It is interesting to note how the curve for \gtol{0.3} sharply increases in slope after time step 40, a moment when the hyperparameters are reoptimized and generate a modified model in which $\gamma_\mrm{tol}$ is crossed more often. This indicates that the $\gamma_\mrm{tol}$ level necessary to achieve a given sampling frequency is also influenced by the hyperparameter values.

\begin{figure}
\centering
\begin{subfigure}[b]{0.45\textwidth}
\includegraphics[scale=0.9]{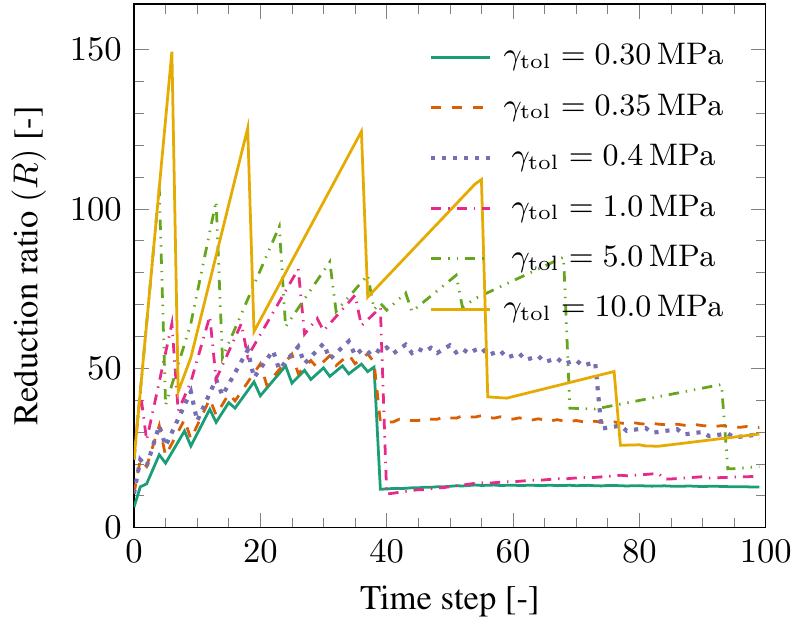}
\end{subfigure}
\begin{subfigure}[b]{0.45\textwidth}
\includegraphics[scale=0.9]{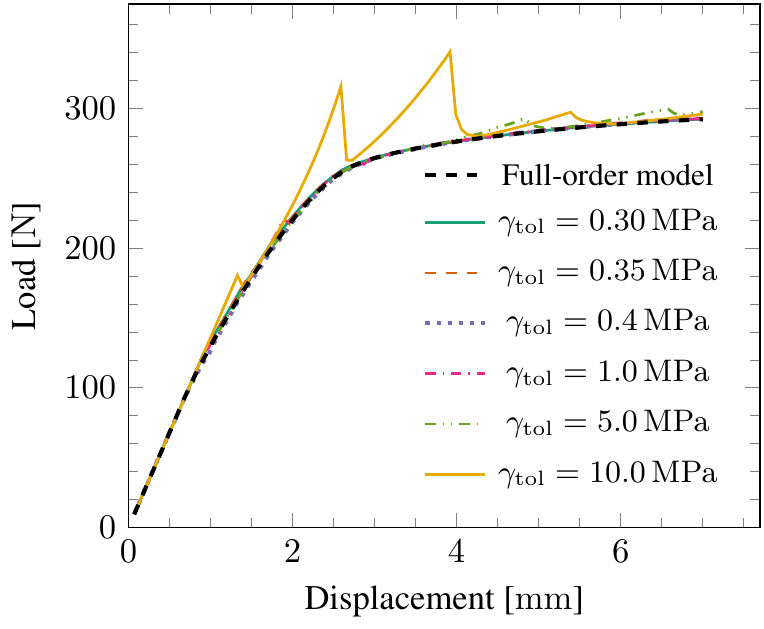}
\end{subfigure}
\caption{Evolution of the reduction ratio and the number of fully-solved steps for models with different values of $\gamma_\mrm{tol}$. Model efficiency cannot be directly controlled by the uncertainty tolerance.}
\label{FIbarratios}
\end{figure}

Similar observations can be made by looking at the degree of control that can be exerted over the acceleration level provided by the framework, which we quantify through the reduction ratio $R$ of \cref{EQreductionratio}. For each model we compute the evolution of the reduction ratio $R$  with the time steps and plot them in \cref{FIbarratios}a. During the first time steps, $\gamma_\mrm{tol}$ has the expected influence on the acceleration level, with ratios as high as 150 being obtained for \gtol{10.0}. When new data is sampled, which requires a number of full-order computations to be performed, the reduction ratio experiences drops that become sharper for higher values of $\gamma_\mrm{tol}$. From time step 40, as the model approaches a perfectly-plastic regime and the sampling frequency increases (see \cref{FIbardsize}), all models converge to reduction ratios between 12 and 30. Adjusting $\gamma_\mrm{tol}$ in order to achieve a desired acceleration level is therefore not possible. Indeed, using higher values lead to lower accuracy --- as can be seen in the load-displacement curves of \cref{FIbarratios}b --- without consistent gains in efficiency.

\subsubsection{Effect of re-estimating the hyperparameters}

\begin{figure}
\centering
\begin{subfigure}[b]{0.45\textwidth}
\includegraphics[scale=0.9]{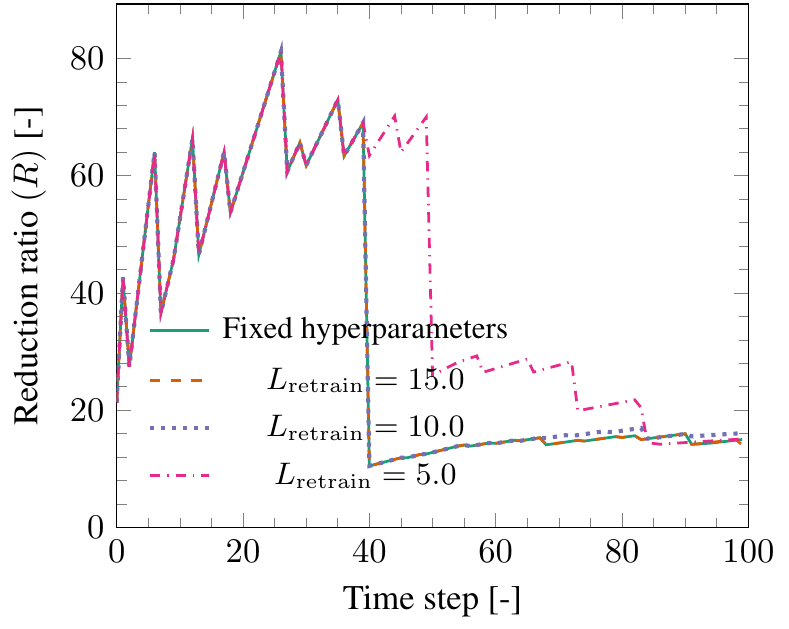}
\caption{\gtol{1.0}}
\end{subfigure}
\begin{subfigure}[b]{0.45\textwidth}
\includegraphics[scale=0.9]{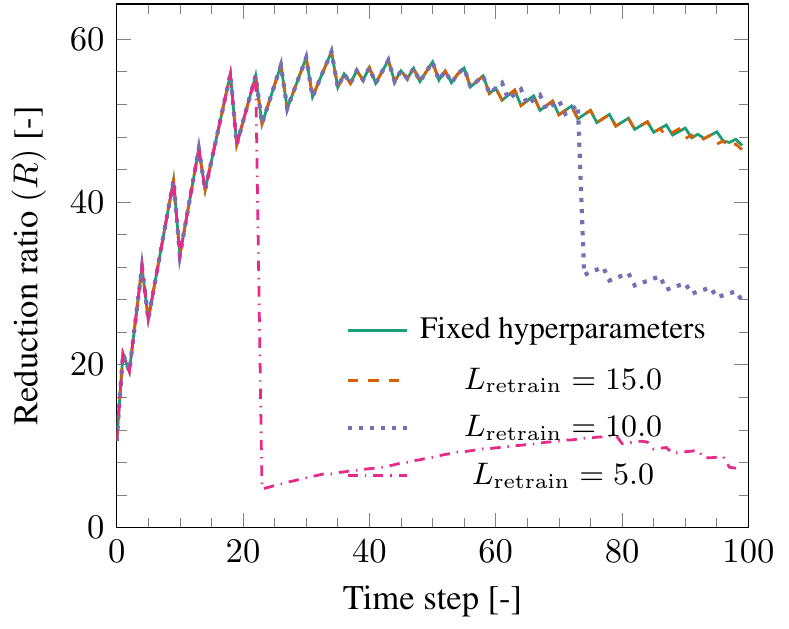}
\caption{\gtol{0.4}}
\end{subfigure}
\caption{Acceleration performance when using the hybrid constitutive model to solve the one-dimensional bar problem for different values of $L_\mrm{retrain}$ and $\gamma_\mrm{tol}$.}
\label{FIbarretrainings}
\end{figure}

Up until this point, the hyperparameters $\noise_\mrm{f}$, $\noise_\mrm{n}$ and $\ell$ have been estimated at the beginning of the analysis and updated with a log marginal likelihood ratio threshold $L_\mrm{retrain}=\SI{10}{}$. For the next example, we return to the one-dimensional model of the previous sections but now using different values of $L_\mrm{retrain}$. Since the sampling frequency dictated by $\gamma_\mrm{tol}$ is also influenced by the hyperparameters, we show results for two different values of $\gamma_\mrm{tol}$. The evolution of the reduction ratio $R$ for these eight models is shown in \cref{FIbarretrainings}.

Here we see two distinct behaviors depending on the uncertainty threshold level. For \gtol{1.0}, allowing for a higher hyperparameter retraining frequency leads to higher acceleration factors up to time step 80. For the lower value of \gtol{0.4}, retraining the hyperparameters has a detrimental effect on efficiency by leading to a much higher sampling frequency, with this change in behavior happening earlier for lower values of $L_\mrm{retrain}$. These results suggest there is no clear recommendation to be made on the optimum hyperparameter retraining frequency, as model performance is dictated by a complex interaction between $\gamma_\mrm{tol}$, the multiple adaptive model components and the hyperparameter values. 

\begin{figure}
\centering
\begin{subfigure}[b]{0.45\textwidth}
\includegraphics[scale=0.9]{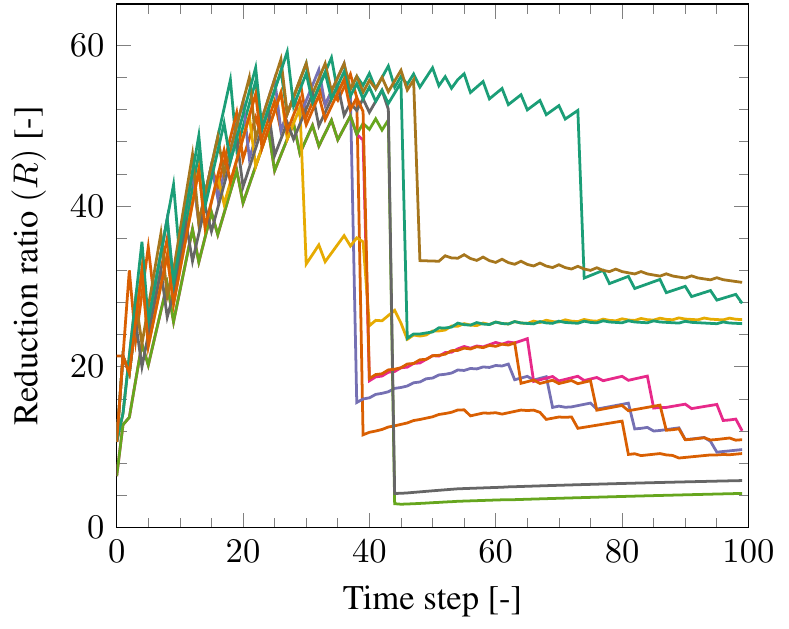}
\end{subfigure}
\begin{subfigure}[b]{0.45\textwidth}
\includegraphics[scale=0.9]{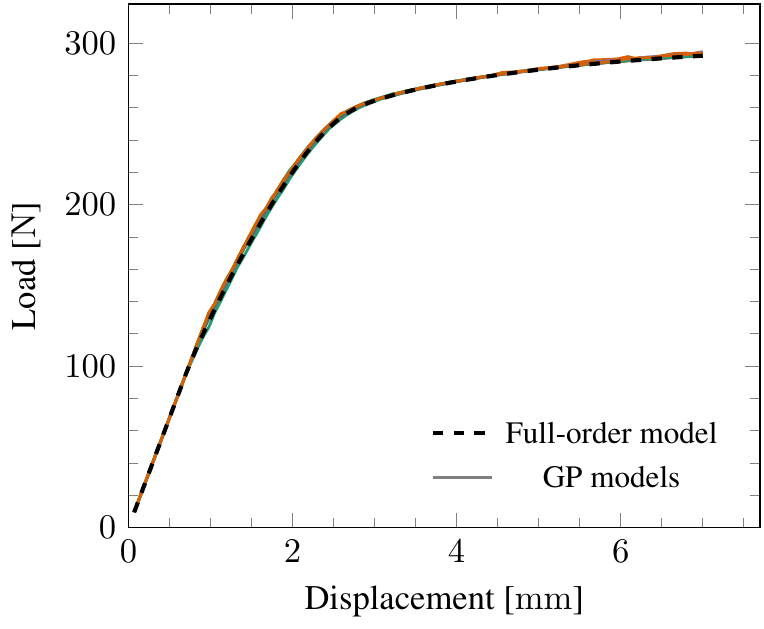}
\end{subfigure}
\caption{Reduction ratios and load-displacement curves for 10 models with different initial hyperparameter values (\gtol{0.4}, $L_\mrm{retrain}=10.0$). A large spread in performance can be observed.}
\label{FIbarhps}
\end{figure}

We show one last example before returning to the two-dimensional version of the model. We now keep both the retraining threshold $L_\mrm{retrain}=10.0$ and the uncertainty tolerance \gtol{0.4} fixed and run 10 models with different seeds being given to the pseudo-random number generator used to initialize the BFGS optimizer that finds the hyperparameter values. Results in terms of reduction ratios and load-displacement curves are shown in \cref{FIbarhps}. With the different initializations, the optimizer finds different local marginal likelihood maxima corresponding to different sets of hyperparameters. This in turn leads to a large spread in acceleration levels between models, once again demonstrating the sensitivity of the framework to the hyperparameter values and to their interaction with $\gamma_\mrm{tol}$. Nevertheless, all 10 models approximate the reference response with excellent accuracy due to the fairly strict value adopted for $\gamma_\mrm{tol}$.

\subsubsection{Acceleration versus mesh density}

We now return to the two-dimensional model shown in \cref{FIbarmesh} in order to investigate how model performance scales with the level of mesh discretization. We solve the model with multiple different mesh densities with characteristic element sizes ranging from \SI{16}{\mm} (32-element mesh) to \SI{0.8}{\mm} (3020-element mesh). All meshes are composed of constant-strain triangles with one integration point each. The model parameters are the same as in the previous examples and the uncertainty threshold is fixed at \gtol{1.0}. Furthermore, and in order to keep the comparison between meshes as consistent as possible, we first use the 3020-element model to estimate the hyperparameters once at the beginning of the study (by loading monotonically along the $k$-means direction), set them as initial values for all other meshes and keep them fixed by adopting a high value for $L_\mrm{retrain}$. This avoids the possibility of GP surrogates from models with different meshes converging to different local maxima of the likelihood function with distinct behaviors (see \cref{FIbarhps}) due to small fluctuations in stress values.

\begin{figure}
\centering
\includegraphics[scale=1.0]{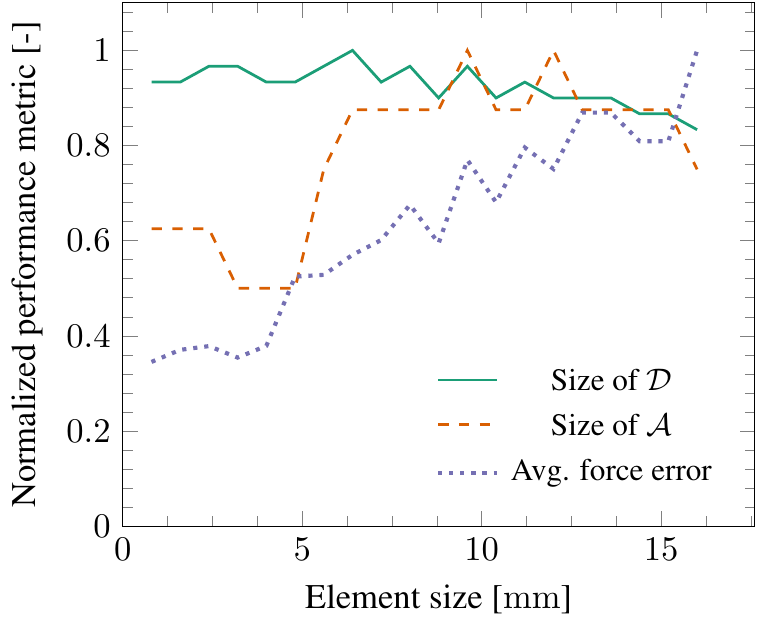}
\caption{Performance of the reduction framework for different levels of mesh discretization (two-dimensional bar problem). The amount of constitutive information necessary to build an accurate surrogate is independent of the total number of integration points in the model.}
\label{FIdogmeshsizemeasures}
\end{figure}

The relative changes in $\vert\dataset\vert$ and $\vert\Ful\vert$ as well as in the average force error with respect to the reference solution are plotted against the discretization level in \cref{FIdogmeshsizemeasures}. It can be seen that the amount of information needed by the model in order to maintain accuracy is independent of the mesh density, with both $\vert\dataset\vert$ and $\vert\Ful\vert$ showing only relatively minor fluctuations as the mesh is refined. This is a consequence of the greedy approach employed here: even though the total number of integration points can be large, points are only sampled or added to \Ful\ if their uncertainty is higher than $\gamma_\mrm{tol}$. 

\begin{figure}
\centering
\includegraphics[scale=1.0]{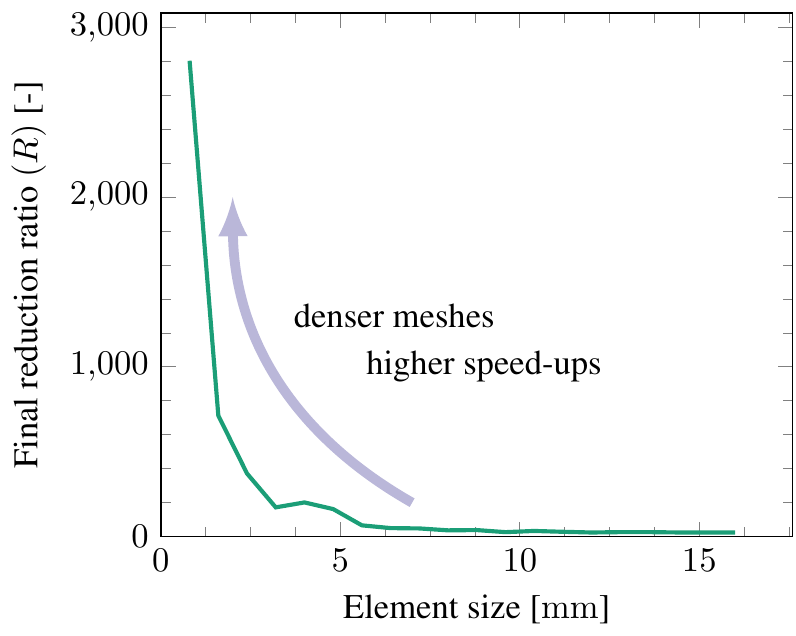}
\caption{Final reduction ratios for the two-dimensional bar problem with different levels of mesh discretization. The efficiency of the reduction framework increases significantly as denser meshes are used.}
\label{FIdogmeshsizeratios}
\end{figure}

Since models with different meshes are approximating the same underlying solution, it is intuitive to expect the sampling effort to be similar: even though the constitutive manifold is more densely evaluated if a denser mesh is used, these evaluations consist of closely-packed clusters in strain space whose response can be accurately approximated by a small number of GP observations\footnote{This is a natural consequence of our choice of kernel, which quantifies the similarity between points by their distance in strain space.}. As a consequence, the reduction ratio $R$ (and therefore the speed-up) increases dramatically with mesh density, as can be seen in \cref{FIdogmeshsizeratios}. Recalling that the same example has been used with \fetwo\ in \cref{SEfe2demonstration}, we can therefore expect that opting for the densest mesh used here would lead to a reduced model almost 3000 times faster than its full-order counterpart without resorting to \emph{offline} training.

\subsubsection{Initial number of fully-solved points}

\begin{figure}
\centering
\begin{subfigure}{\textwidth}
\centering
\includegraphics[scale=1]{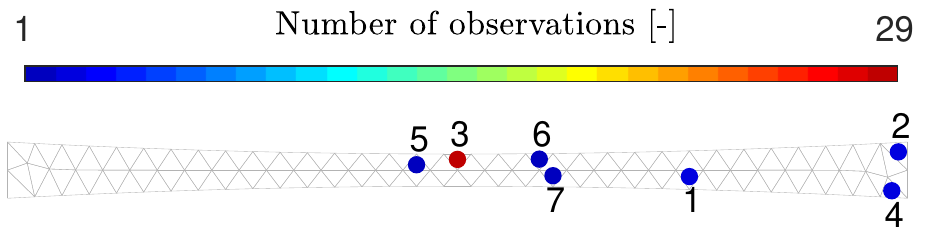}
\caption{$k=1$}
\end{subfigure}
\begin{subfigure}{\textwidth}
\centering
\includegraphics[scale=1]{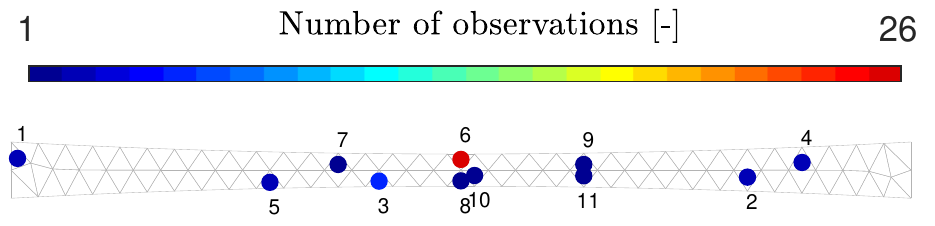}
\caption{$k=5$}
\end{subfigure}
\begin{subfigure}{\textwidth}
\centering
\includegraphics[scale=1]{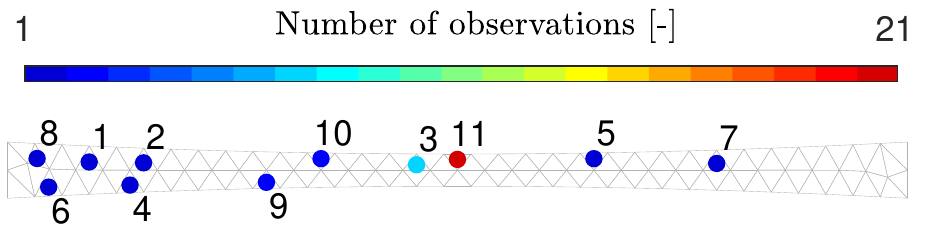}
\caption{$k=10$}
\end{subfigure}
\caption{Heatmaps of the total number of samplings of points in \Ful\ for different values of $k$. The numbers next to the points indicate the order at which the points are added to \Ful. For each case, the first $k$ points are initially present.}
\label{FIdogclusterheatmaps}
\end{figure}

As one final parametric study on the two-dimensional model of \cref{FIbarmesh}, we investigate the effect of changing the clustering parameter $k$ that determines the initial number of points in \Ful. \cref{FIdogclusterheatmaps} shows a set of heatmaps plotting the total number of times each anchor point is sampled during the analysis for three different values of $k$. For $k=1$, we start with a point midway between the load application face and the center of the bar and the greedy data selection approach promptly locates the point at the center of the bar undergoing the largest strains. Additional points are eventually added due to a number of canceled steps, but the model nevertheless concentrates most of the sampling effort at the center and the remaining points remain dormant for the rest of the analysis. The same happens for the models with $k=5$ and $k=10$: the model starts with $k$ anchor models that remain dormant and immediately adds another one at the center of the bar from where most of the training data is obtained. Due to the relatively simple strain path of this specific example, increasing $k$ does not seem to be beneficial. However, the greedy framework is able to naturally disregard the redundant information and concentrate the sampling effort where it is needed. All three models have, therefore, similar accuracy and acceleration levels.

\subsection{Two-dimensional plate with multiple cutouts}

\begin{figure}
\centering
\includegraphics[scale=0.5]{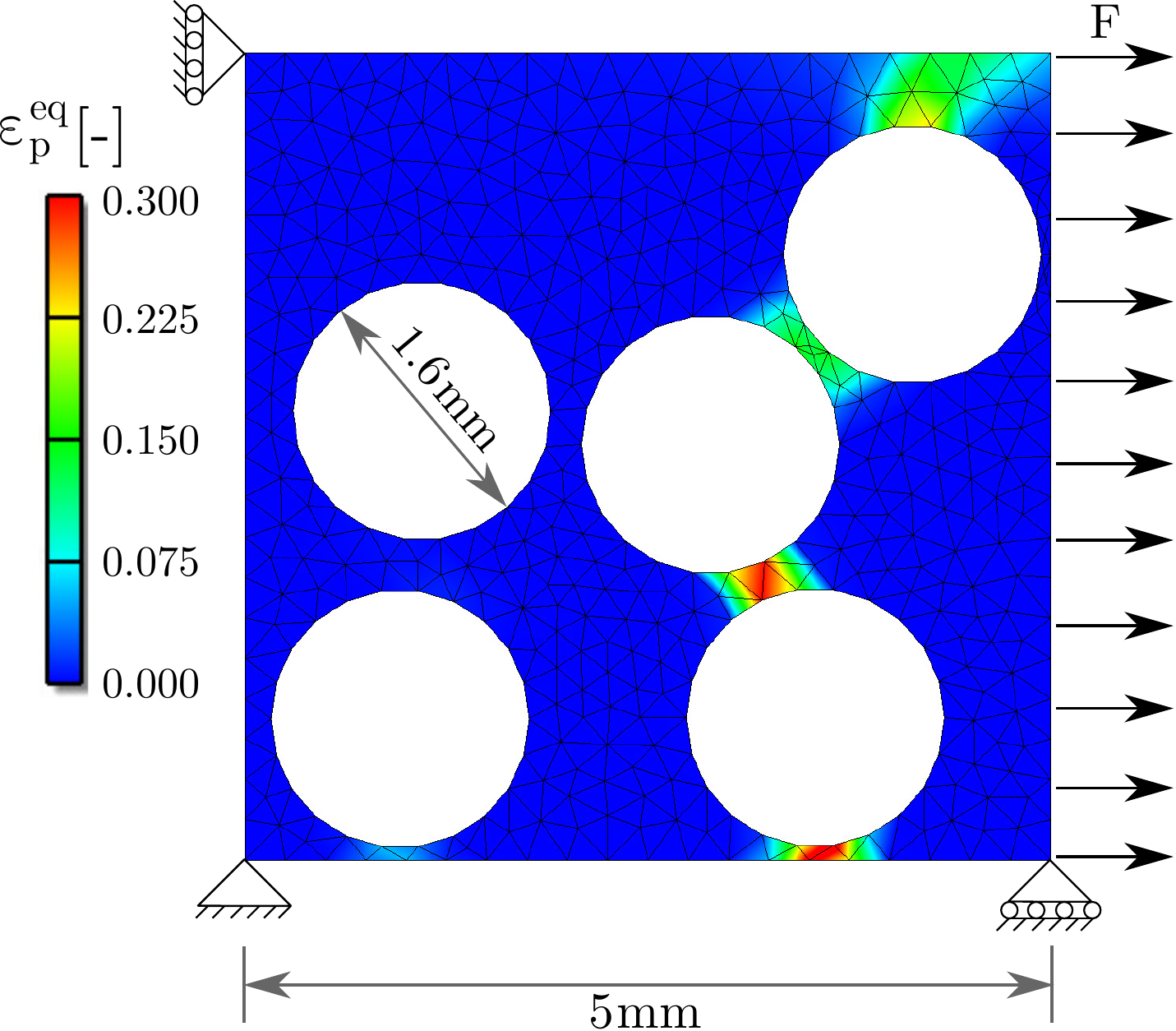}
\caption{Geometry, boundary conditions and final equivalent plastic strain distribution for the two-dimensional cutout example.}
\label{FImultiholeepspeq}
\end{figure}

We now move to an example with complex geometry in order to investigate how the reduction framework fares in approximating a larger portion of the original full-order constitutive manifold. The example employs the same elastoplastic material as before but now concerns the plate with multiple cutouts with boundary conditions and final plastic strain distribution shown on \cref{FImultiholeepspeq}. Due to the presence of the cutouts, the stress distribution is considerably more complex than for the previous examples. As the load increases, plastic strain arises at the stress concentration regions between cutouts and forms a strain localization band spanning the complete height of the model. The plate is discretized with $754$ constant-strain triangles with one integration point each.

\begin{figure}
\centering
\includegraphics[scale=1.0]{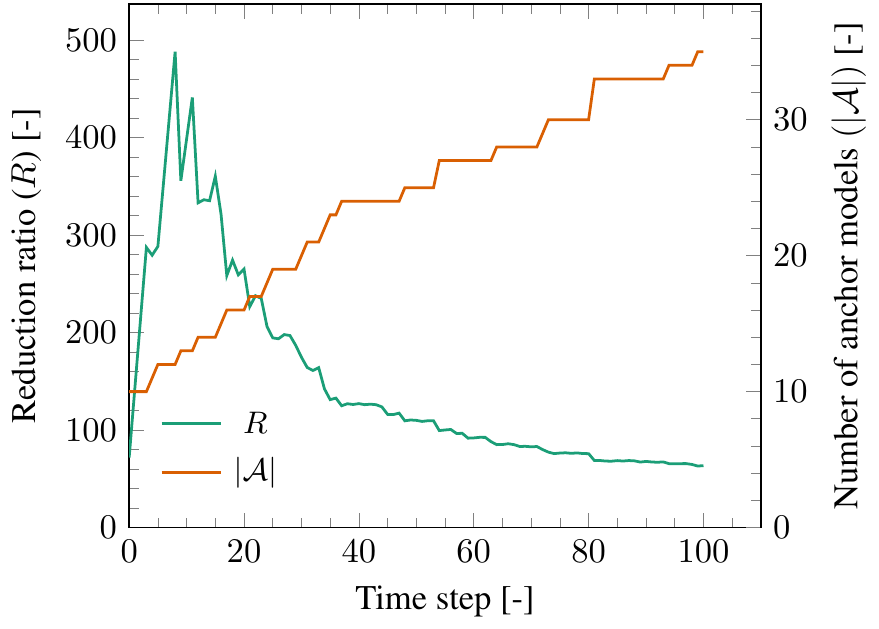}
\caption{Evolution of the reduced-order solution of the two-dimensional plate problem with cutouts. With a more complex constitutive behavior to approximate, the size of the \Ful\ set is larger for this example.}
\label{FImultiholemetrics}
\end{figure}

We solve the problem with $k=\SI{10}{}$, \gtol{2.0}, $L_\mrm{retrain}=\SI{10}{}$ and $\gamma_\mrm{cancel}=\SI{80}{\mega\pascal}$ for a total of 100 time steps. As in the previous examples, the initial hyperparameters are estimated by sampling the initial $k$-means directions. We plot the evolution of the reduction ratio $R$ and the size of the fully-solved set \Ful\ in \cref{FImultiholemetrics}. While the model of \cref{FIdogclusterheatmaps} is able to rely on the information coming from a single anchor model to accurately describe its constrained constitutive space, the complex stress distribution of the current example demands the sampling of a significantly higher number of points. The acceleration is therefore smaller than for the previous cases, with a reduction ratio of approximately \SI{63}{} at the end of the analysis. This result is not unexpected, since the reduction framework relies on the assumption that the macroscopic geometry and boundary condidions constrain the constitutive response to lie on a manifold of much lower complexity. The complex stress state treated here challenges this assumption. 

\begin{figure}
\centering
\begin{subfigure}[b]{0.45\textwidth}
\includegraphics[scale=0.7]{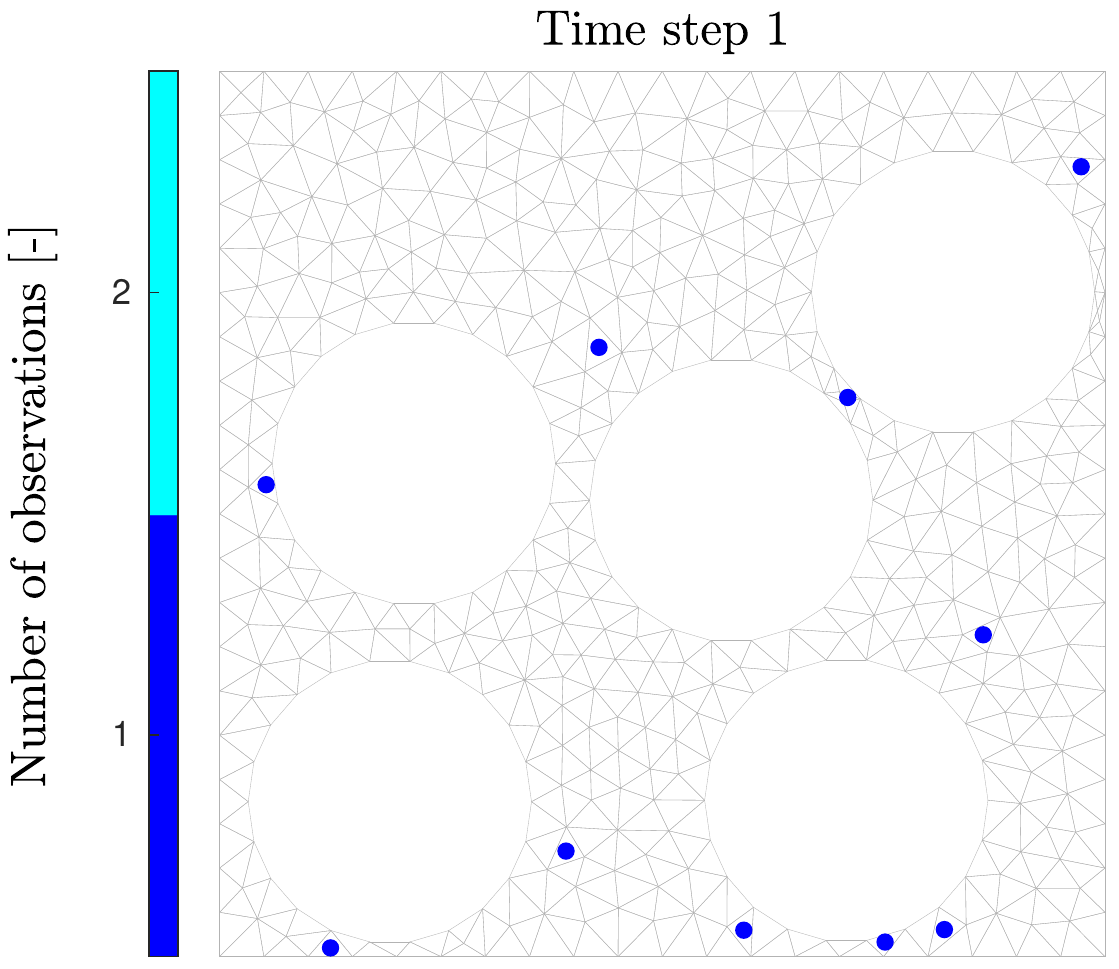}
\end{subfigure}\hfill
\begin{subfigure}[b]{0.45\textwidth}
\includegraphics[scale=0.7]{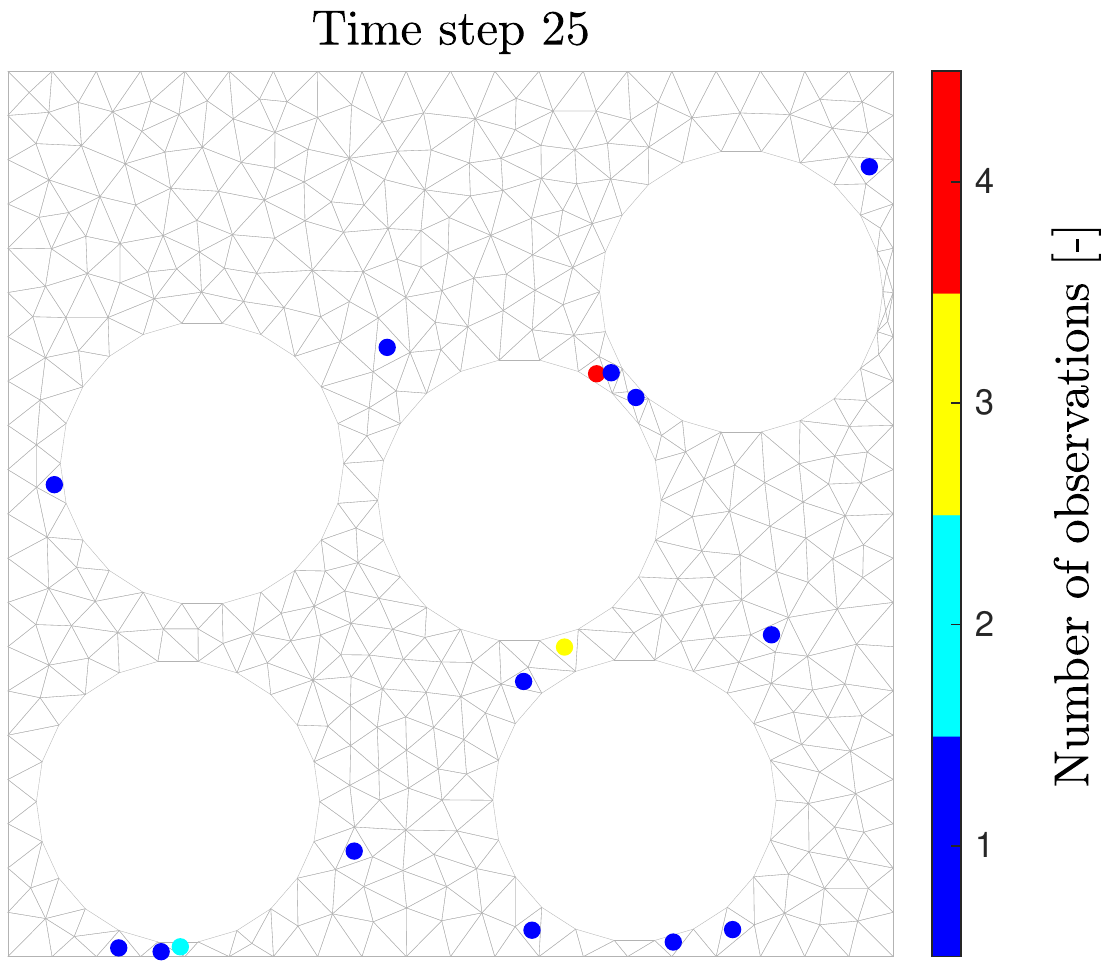}
\end{subfigure}
\begin{subfigure}[b]{0.45\textwidth}
\includegraphics[scale=0.7]{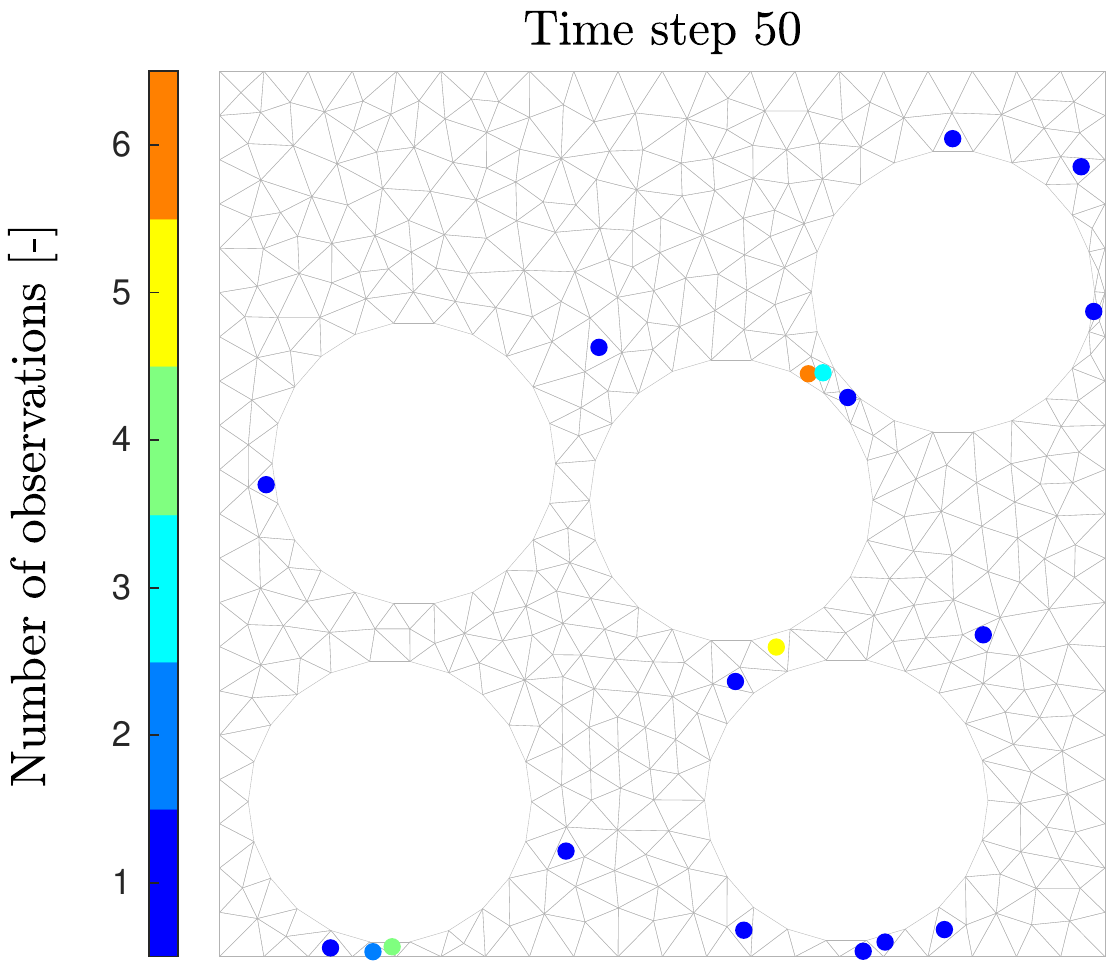}
\end{subfigure}\hfill
\begin{subfigure}[b]{0.45\textwidth}
\includegraphics[scale=0.7]{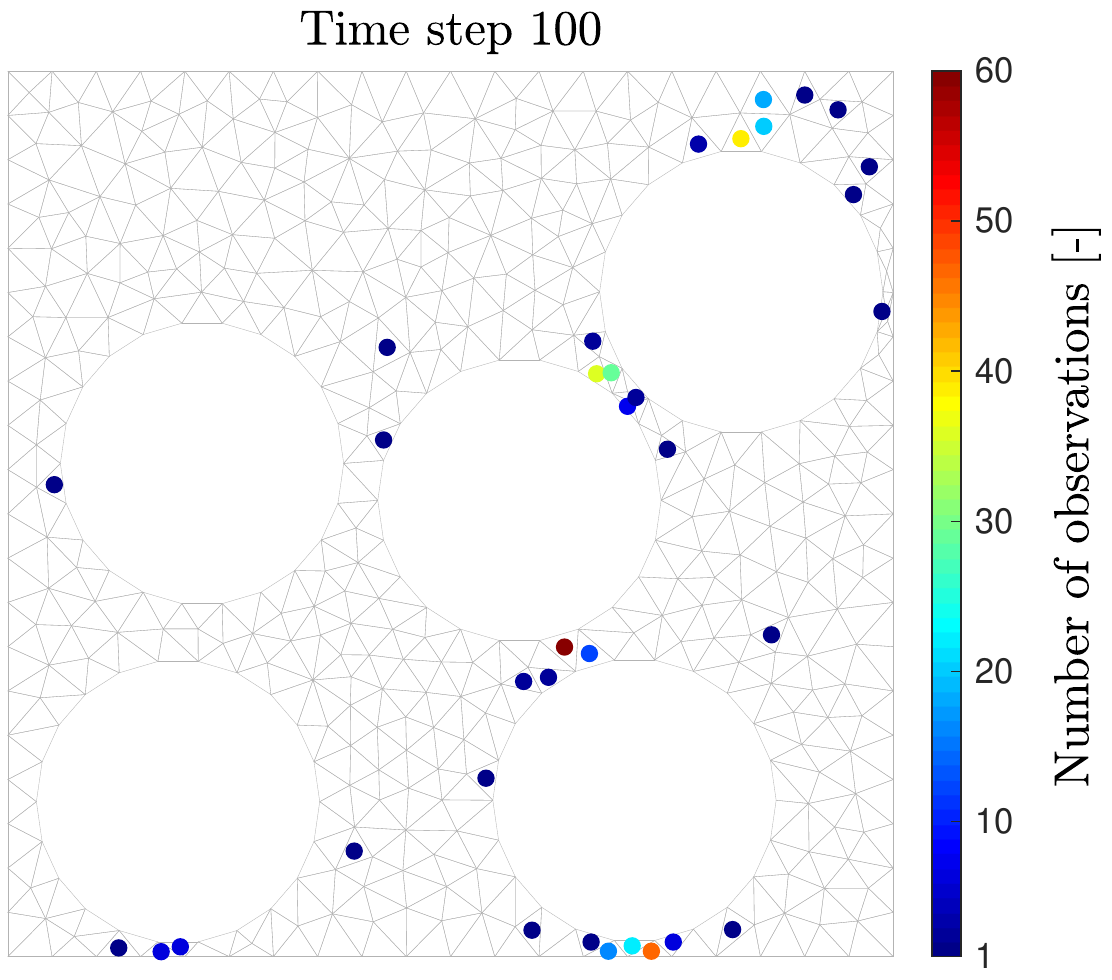}
\end{subfigure}
\caption{Heatmaps of the total number of samplings of points in \Ful\ for the two-dimensional cutout example taken at four different moments throughout the analysis. The adaptive model automatically concentrates the sampling effort at regions undergoing strain localization.}
\label{FImultiholeheatmap}
\end{figure}

\begin{figure}
\centering
\includegraphics[scale=1.0]{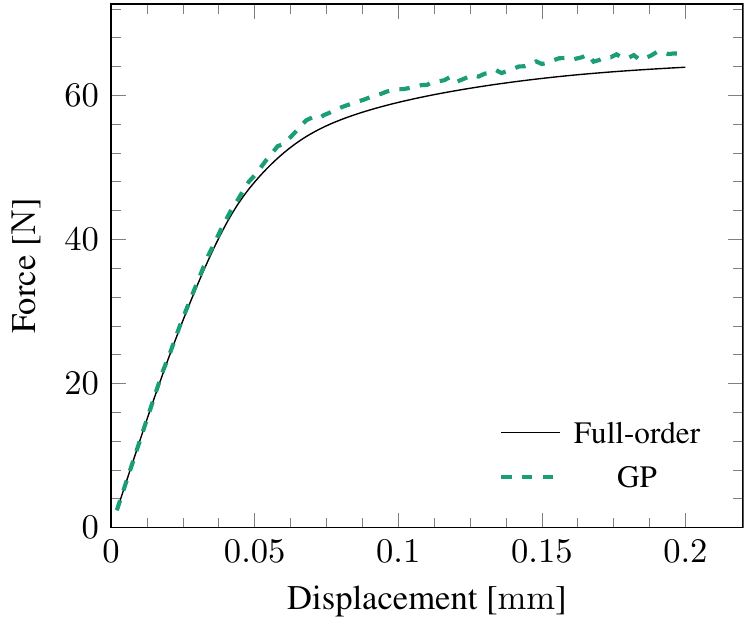}
\caption{Load-displacement curves for the two-dimensional plate with cutouts obtained with the full- and reduced-order models.}
\label{FImultiholelodi}
\end{figure}

It is also interesting to plot the heatmap of GP observations in order to visualize which points are being sampled. Results can be seen in \cref{FImultiholeheatmap}. In contrast with the model of \cref{FIdogclusterheatmaps}, sampling is performed on a larger number of points. This indicates that different parts of the mesh experience significantly different strain paths. Comparing \cref{FImultiholeheatmap} with the plastic strain field of \cref{FImultiholeepspeq}, we see that the point distribution is closely related to how plastic strain is distributed throughout the domain. The framework is therefore able to direct computational effort to regions in the mesh where it is most needed while employing an efficient approximation for the rest of the domain. Finally, the load-displacement curves of the full-order and reduced models are shown in \cref{FImultiholelodi}. Again a satisfactory agreement is obtained at a fraction of the number of full model evaluations.

It is worth mentioning, however, that this more complex example is less numerically stable than the previous ones. Running the model with higher values of $\gamma_\mrm{tol}$ leads to convergence issues as new data coming from the anchor models leads to large jumps between equilibrium solutions that push the stability of the Newton-Raphson solver to its limit. Additionally, a model with reasonable stability could only be obtained by forcibly increasing the smoothness of the GP approximation by increasing the lower bound of $\sigma^2_\mrm{n}$ to \SI{1.0}{\mega\pascal\squared}. Further research effort is therefore necessary in order to improve the stability of the present approach when faced with highly-complex strain distributions.

\subsection{Mixed-mode cohesive crack propagation}

We close the present discussion with one final example exploring the use of the framework to approximate a traction-separation response. In \cref{SEcohesivehomogenization}, we argue that bulk homogenization in \fetwo\ loses objectivity upon global softening at the microscale, at which point switching to a cohesive homogenization strategy becomes necessary. Here we construct a surrogate model for the associated macroscopic cohesive material by using the framework of \cref{SEsurrogatemodelingframework} but now defining $\mathcal{S}$ as:
\begin{equation}
\bs{\tau}^\Omega\left(\bs{\tau}^\Omega_\mrm{eff},\dataset\right) = \mathbb{E}\left[\widehat{\bs{\tau}}\left(\bs{\tau}^\Omega_\mrm{eff},\dataset\right)\right]
\end{equation}
\noindent where $\bs{\tau}^\Omega_\mrm{eff}$ is an effective traction computed from the displacement jump $\mbf{\jump}$ and the stress at crack initiation that accounts for the singular nature of the initially-rigid cohesive law and takes the place of the shifted jump ${\llbracket \mbf{v}\rrbracket}^\Omega$ of \cref{EQfe2traction} as input variable \cite{vdmeer09ijf}. Note that, in contrast to \cref{EQcorrectionmodel}, assuming a correction from elasticity ceases to be interesting here and we therefore build a direct regression for the tractions with one GP model for each component of the surrogate traction vector $\bs{\tau}^\Omega$. Furthermore, we exploit the knowledge that decohesion is an irreversible process by switching to the trivial solution $\bs\tau=\mbf{0}$ after the GP confidently predicts zero traction for the first time:
\begin{equation}
\text{for a given point}\,p\text{, if}\,\norm{\mathbb{E}\left[\widehat{\bs{\tau}}\right]}_p\leq\sigma_\mrm{n}\,\,\mrm{and}\,\,\gamma^\mrm{o}_p\leq\gamma_\mrm{tol}\,\,\Rightarrow\,\,\text{switch to }\,\,\bs{\tau}^\Omega=\mbf{0}\,\vert\,\forall\,\bs{\tau}^\Omega_\mrm{eff}>\mbf{0}
\label{EQgpzerotraction}
\end{equation}
\noindent After the switch, we stop using the GP approximation for the integration point in question. This modification improves efficiency because we avoid having to train the GP to reproduce the fully-damaged branch of the cohesive law\footnote{The GP prediction naturally moves to its zero-traction prior away from \dataset, but retraining would still be periodically triggered due to the variance increase.}.

\begin{figure}
\centering
\includegraphics[width=0.7\textwidth]{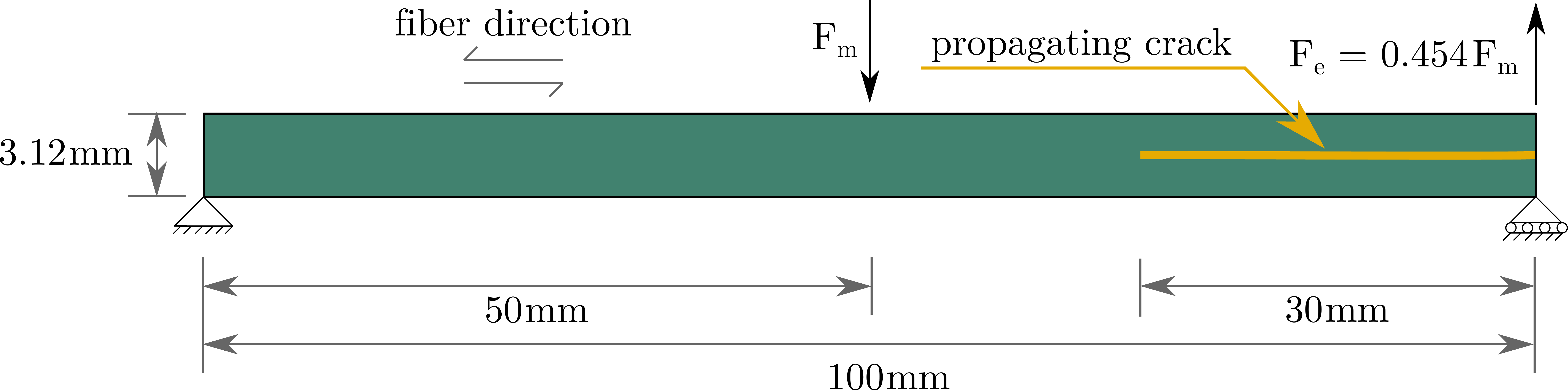}
\caption{Mixed-mode crack propagation example: Geometry, loads and boundary conditions.}
\label{FImmbmodel}
\end{figure}

The example concerns the mixed-mode bending test shown in \cref{FImmbmodel} and is taken from \cite{vdmeer09ijf}, where the same structure is solved for multiple mode-mixity ratios. Here we opt for a single ratio $\alpha=G_\mrm{II}/\left(G_\mrm{I}+G_\mrm{II}\right)=0.5$ and therefore only deal with the specific ratio between applied forces shown in \cref{FImmbmodel} and with a single initial notch length of \SI{30}{\mm}. In an \fetwo\ approach, both the bulk material behavior and the cohesive softening response would be derived from embedded RVEs by employing \cref{EQfe2stress,EQfe2traction}. Without loss of generality and in keeping with the original model of \cite{vdmeer09ijf}, in this demonstration we instead use a linear-elastic orthotropic model for the bulk material and a bulk stress-constrained cohesive zone law to model the traction-separation behavior of the propagating crack. The model is solved in plane stress and initially discretized with 3107 4-node quadrilateral elements with 4 integration points each. Cohesive segments are inserted on the fly by using the Phantom Node method \cite{Hansbo2004} (later renamed CutFEM \cite{Burman2015}), with elements being duplicated in order to describe a displacement jump running through the elements as the crack propagates from the tip of the notch.

\begin{figure}
\centering
\includegraphics[scale=1.0]{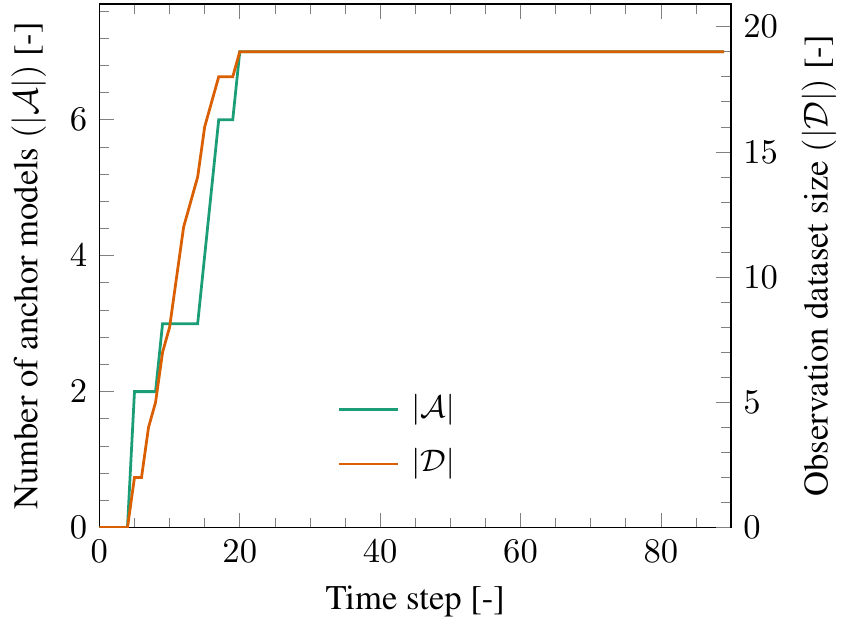}
\caption{Mixed-mode crack propagation example: Evolution of the learning process in terms of the sizes of the \Ful\ and \dataset\ sets. No learning occurs during most of the crack propagation process. The framework also refrains from computing the full-order cohesive model for most of the analysis.}
\label{FImmbmetrics}
\end{figure}

In order to keep the discussion simple, we only use the GP framework to approximate the response at cohesive integration points, but extending the example to also use GP for the bulk response would be straightforward. We run the reduced model with \gtol{1.0}, $\gamma_\mrm{cancel}=\SI{100}{\mega\pascal}$ and fixed hyperparameters obtained from sampling a fictitious anchor model in the initial effective traction direction seen by the first cohesive point. Note that in this case we cannot rely on the $k$-means strategy from the previous examples since the analysis starts with no cohesive integration points.

\begin{figure}
\centering
\includegraphics[scale=1.0]{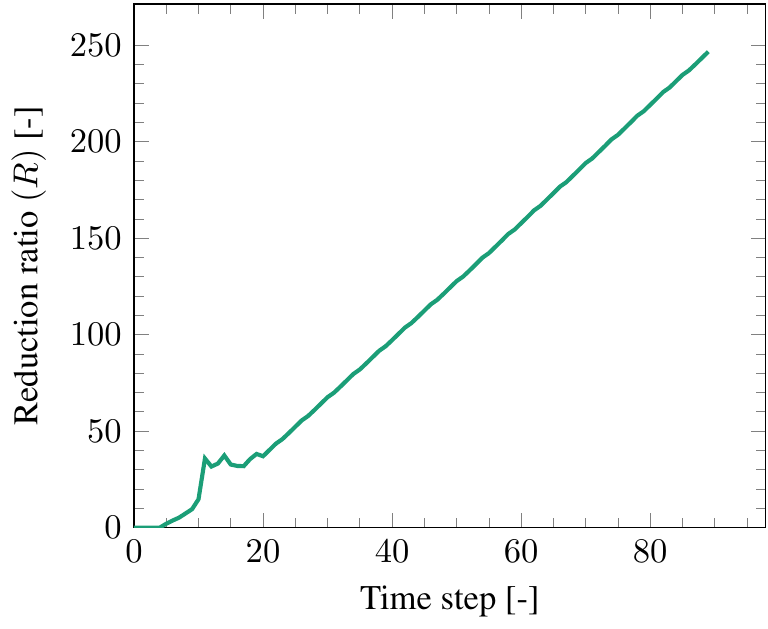}
\caption{Evolution of the reduction ratio $R$ for the mixed-mode crack propagation example. The fact that \dataset\ is not growing for most of the analysis allows for high acceleration ratios to be obtained.}
\label{FImmbratio}
\end{figure}

The solution process is tracked by plotting the evolution of the \Ful\ and \dataset\ sets on \cref{FImmbmetrics}. During the first five time steps, the structure is loaded until the onset of crack propagation. As no cohesive points exist at this stage, \Ful\ and \dataset\ remain empty. The first cohesive integration points created when the crack starts to propagate are added to \Ful\ and their responses are sampled into \dataset. From that moment on, the framework is left to decide which points are added. Since the crack tip can be seen as a moving source travelling through the domain, points created after time step 20 can be accurately approximated with information obtained from points closer to the notch. Since the decohesion process follows almost exactly the same path in these subsequent integration points and we switch to a trivial solution with zero traction for fully-damaged points (\cref{EQgpzerotraction}), no new data is needed in \dataset\ for the rest of the analysis. The greedy algorithm is able to detect this and interrupt the learning process, with the 7 anchor models in \Ful\ remaining dormant for the rest of the analysis and no more full-order updates being computed. This results in the high values for the reduction ratio plotted in \cref{FImmbratio}.

\begin{figure}
\centering
\includegraphics[scale=1.0]{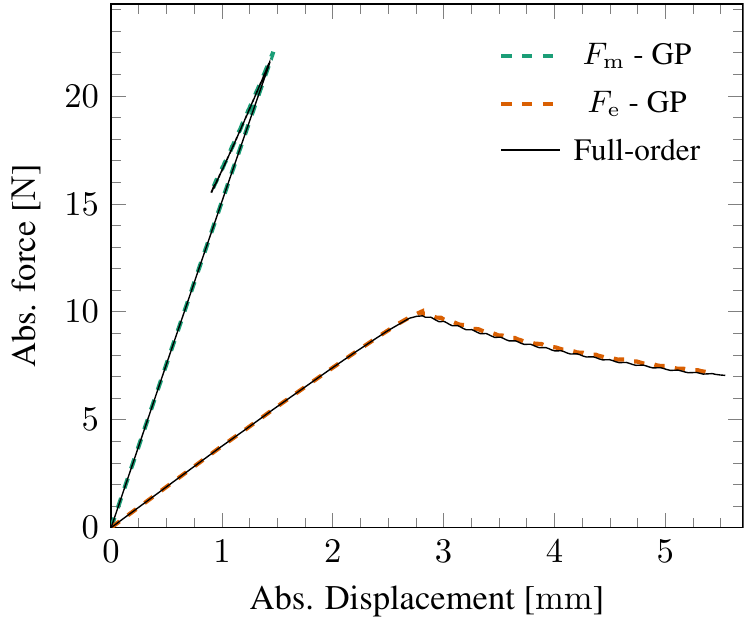}
\caption{Load-displacement curves for the mixed-mode crack propagation example. The active learning framework shows excellent agreement with the full-order response.}
\label{FImmblodi}
\end{figure}

\begin{figure}
\centering
\includegraphics[scale=1.0]{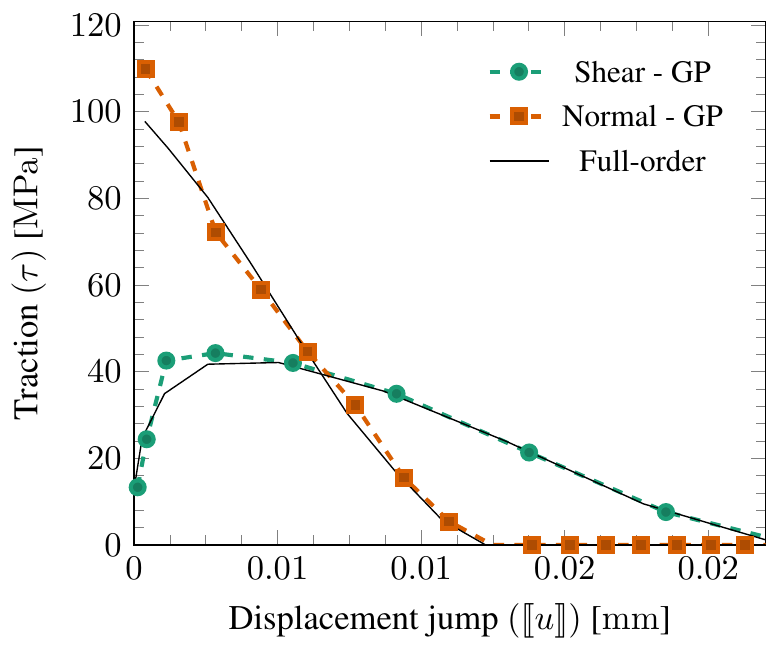}
\caption{Evolution of the mixed-mode decohesion of an integration point computed by the GP model during the analysis phase on which $\dataset$ is not growing.}
\label{FImmbtractions}
\end{figure}

The load-displacement curves obtained with the full and hybrid models are shown in \cref{FImmblodi}, where we plot absolute forces and displacements from the two load locations shown in \cref{FImmbmodel}. It can be seen that the global response obtained with the active learning approach is virtually indistinguishable from the full-order one. Finally, we can observe how accurate the local traction approximation is by plotting in \cref{FImmbtractions} the traction-separation curves for a cohesive integration point created and completely solved after the anchor models become dormant. We see that the complete mixed-mode softening process is correctly predicted by the GP models based solely on the response of earlier points.

\section{Conclusions}
\label{SEconclusions}

This work introduces an adaptive probabilistic learning framework for the \emph{online} construction of surrogate constitutive models for concurrent multiscale analysis. The framework eliminates the need to sample a potentially infinitely-dimensional input space \emph{offline} and instead fits a set of Gaussian Process (GP) models with data sampled \emph{online} from a small number of fully-solved \emph{anchor} models. The approach incorporates additional physics information by enhancing the conventional GP regression with tangent stiffness observations coming from the anchor models. A greedy data selection procedure ensures the surrogate response is kept accurate, efficient and independent of the time step size and of the macroscopic discretization level.

The reduction approach was described in detail and its performance was assessed with an extensive set of numerical examples. The ability of the framework of reducing the computational effort associated with \fetwo\ was demonstrated with a preliminary example, after which a detailed parametric study was performed on single-scale models without loss of generality. The uncertainty tolerance parameter used to control the GP sampling frequency was found to provide only a limited degree of control on the balance between accuracy and efficiency of the reduced response due to the presence of other model components that can also trigger a refinement of the GP approximations (\eg canceled solutions). Using a larger initial set of fully-solved points or allowing for the GP hyperparameters to be re-estimated during the analysis were found to exert little influence on the performance of the model, at least for the specific examples treated here. The acceleration brought by the greedy learning strategy was found to drastically increase as the macroscopic mesh is refined, with reduction ratios as high as 2800 times being obtained.

An additional example was used to demonstrate the ability of the reduction approach to handle models with complex stress distributions. Although the acceleration was lower in this case due to the complexity of the constitutive manifold being approximated, the greedy sampling strategy successfully concentrated the learning effort on the most informative mesh regions. One final model involving mixed-mode crack propagation was presented. In contrast with the previous examples, the nature of the crack propagation problem allowed the GP to reuse previously obtained information rather than continuously growing the dataset. Acceleration ratios of up to 250 times were obtained and the GP surrogate was able to take over the entire set of integration points from the original full-order model for most of the analysis. 

The presented results suggest the framework is a promising approach in reducing the computational effort of nonlinear concurrent multiscale modeling and circumventing the curse of dimensionality associated with the \emph{offline} construction of surrogates for path-dependent materials. Nevertheless, further model development is necessary in order to allow for non-monotonic load paths including unloading/reloading behavior since the adopted GP formulation assumes a unique mapping between strains and stresses. Although path dependency is accounted for in the sense that macroscopic models with different strain paths will construct different surrogates, an accurate approximation is not guaranteed for models that visit points in strain space more than once during their loading paths.

\section*{Acknowledgements}

The authors gratefully acknowledge financial support from the Netherlands Organization for Scientific Research (NWO) under Vidi grant nr. 16464. 

\bibliographystyle{elsarticle-num}
\bibliography{paper}

\end{document}